\magnification 1100

  
  %
  %

  %
  %
  \let\footnotea=\footnote
  \def\anote#1#2{\footnotea{\hbox{$^{#1}$}}{\eightpoint#2}}  
  \catcode`@=12 

 \def\defrefnote#1{\definexref{#1}{{\the\footnotenumber}}{refnotes}}

  %
  %


 \input eplain.tex
\makeatletter
\def\numberedfootnote{%
ÊÊ\global\advance\footnotenumber by 1
ÊÊ\@eplainfootnote{{\number\footnotenumber}}%
}%
\def\makecolumns#1/#2 {\par \begingroup
ÊÊ \@columndepth = #1
ÊÊ \advance\@columndepth by -1
ÊÊ \divide \@columndepth by #2
ÊÊ \advance\@columndepth by 1
ÊÊ \@linestogoincolumn = \@columndepth
ÊÊ \@linestogo = #1
ÊÊ \currentcolumn = 1
ÊÊ \def\@endcolumnactions{%
ÊÊÊÊÊÊ\ifnum \@linestogo<2
ÊÊÊÊÊÊÊÊ \the\crtok \egroup \endgroup \par 
ÊÊÊÊÊÊ\else
ÊÊÊÊÊÊÊÊ \global\advance\@linestogo by -1
ÊÊÊÊÊÊÊÊ \ifnum\@linestogoincolumn<2
ÊÊÊÊÊÊÊÊÊÊÊÊ\global\advance\currentcolumn by 1
ÊÊÊÊÊÊÊÊÊÊÊÊ\global\@linestogoincolumn = \@columndepth
ÊÊÊÊÊÊÊÊÊÊÊÊ\the\crtok
ÊÊÊÊÊÊÊÊ \else
ÊÊÊÊÊÊÊÊÊÊÊÊ&\global\advance\@linestogoincolumn by -1
ÊÊÊÊÊÊÊÊ \fi
ÊÊÊÊÊÊ\fi
ÊÊ }%
ÊÊ \makeactive\^^M
ÊÊ \letreturn \@endcolumnactions
ÊÊ \@columnwidth = \hsize
ÊÊÊÊ \advance\@columnwidth by -\parindent
ÊÊÊÊ \divide\@columnwidth by #2
ÊÊ \penalty\abovecolumnspenalty
ÊÊ \noindent 
ÊÊ \valign\bgroup
ÊÊÊÊ &\hbox to \@columnwidth{\strut \hsize = \@columnwidth ##\hfil}\cr
}%
\makeatother

\lefteqnumbers
   \def\testd{oui}
   \def\choixlat{\ifx\numadroite\testd\righteqnumbers
            \else  \lefteqnumbers\fi}
    \choixlat

\catcode`@=\letter
\def\@eplainfootnote#1{\let\@sf\empty 
  \ifhmode\edef\@sf{\spacefactor\the\spacefactor}\/\fi
  \global\advance\hlfootlabelnumber by 1
  \hlstart@impl{foot}{\hlfootlabel}%
  \hldest@impl{footback}{\hlfootbacklabel}%
  \hbox{$^{(#1)}$}%
  \hlend@impl{foot}%
  \@sf\vfootnote{#1.}%
}%
\catcode`@=\other

  \interfootnoteskip=0pt
  \let\note=\numberedfootnote
  \everyfootnote={\eightpoint\leftskip=5truemm\rightskip5truemm}
  
  \hsize150truemm\vsize 240truemm\hoffset=5truemm
  \def\dimstand{\hsize 150truemm\vsize 240truemm\hoffset=5truemm\voffset=0truemm}
  
  \pretolerance=500\tolerance=1000\brokenpenalty=5000
  \parindent3mm
  
  \countdef\temps=170
  \temps=\time
  \countdef\nminutes=171{\nminutes = \time}
  \countdef\nheure=172
  \def\heure{\begingroup                     
     \temps = \time \divide\temps by 60
     \nheure = \temps                        
     \nminutes = \time
     \multiply\temps by 60
     \advance\nminutes by -\temps            
     \ifnum\nminutes<10 \toks1 = {0}%
     \else\toks1 = {}%
     \fi
     \number\nheure h\the\toks1 \number\nminutes  
  \endgroup}%

  \newcount\chstart
  \chstart=\pageno
 \headline={\ifnum\pageno=\chstart {\hfill} \else {\hss \tenrm --\ \folio\ --\hss}\fi}
  \footline={\hfill}
  \normalbaselines
  \frenchspacing
    \def\dater{\vglue-10mm\rightline{(\the\day/\the\month/\the\year)}}
  \def\dateheure{\vglue-10mm\rightline{(\the\day/\the\month/\the\year,\ \heure)}}

  \newif\ifpagetitre \pagetitretrue
\newtoks\hautpagetitre \hautpagetitre={\hfill}
\newtoks\baspagetitre \baspagetitre={\hfill}
\newtoks\auteurcourant \auteurcourant={\hfill}
\newtoks\titrecourant \titrecourant={\hfill}
\newtoks\hautpagegauche
\newtoks\hautpagedroite
\newtoks\hautpagemilieu
\hautpagemilieu={\tenrm\hfil -- \folio\ -- \hfil}
\hautpagegauche={\ifx\midfolio\oui\the\hautpagemilieu\else\tenrm\folio\hfill\the\auteurcourant\hfill\fi}
\hautpagedroite={\ifx\midfolio\oui\the\hautpagemilieu\else\hfill\the\titrecourant\hfill\tenrm\folio\fi}
\newtoks\baspagegauche \baspagegauche={\hfil}
\newtoks\baspagedroite \baspagedroite={\hfil}
\headline={\ifpagetitre\the\hautpagetitre
\else\ifodd\pageno\the\hautpagedroite\else\the\hautpagegauche\fi\fi }
\footline={\ifpagetitre\the\baspagetitre
\else\ifodd\pageno\the\baspagedroite
\else\the\baspagegauche\fi\fi \global\pagetitrefalse}

\def\pageblanche{\vfill\eject\pagetitretrue
\null\vfill\eject
\pagetitretrue
}
\def\chgtpage{\ifodd\pageno \else
\pageblanche \fi \pagetitretrue\titreun=0\footnotenumber=0}

\def\chgtpageincrtitreun{\ifodd\pageno \else
\pageblanche \fi \pagetitretrue\footnotenumber=0}

\def\majnombres{\ifodd\pageno \else
\pageblanche \fi \pagetitretrue\hautpoly\titreun=0\footnotenumber=0}

\def\hautspages#1#2{\auteurcourant={\ninepcap#1}\titrecourant={\nineit#2}}

\ifnum\chstart=\pageno \pagetitretrue\fi
  


  \def\leftnote#1{\vadjust{\setbox1=\vtop{\hsize 20mm\parindent=0pt\eightpoint
  \baselineskip=9pt\rightskip=4mm plus 4mm\vskip-4.7mm#1}\hbox{\kern-2cm\smash{\box1}}}}

  
  \def\raggedcenter{\leftskip=20pt plus 10em  
       \rightskip=\leftskip 
        \parfillskip=0pt 
         \spaceskip=.3333em \xspaceskip=.5em 
          \pretolerance=9999 \tolerance=9999
           \hyphenpenalty=9999 \exhyphenpenalty=9999 }
           
  \def\titrecentre#1{{\parindent0mm\raggedcenter
       \spaceskip=.6em plus .2em minus .2em\xspaceskip=.6em plus .2em minus .2em
        \tit#1\par}}
        


  \def\oui{oui}
  
   \def\fontetitreun{\ifx\paradouze\oui\douzepts\gpdouze\twelvebf\textfont1=\twelveib\else
\quatorzepts\gpquatorze\fourteenbf\fi}

\def\fontetitreunl{\douzepts\textfont1=\twelveib\scriptfont1=\tenib\fourteenti}
 
 \def\fontetitredeux{\textfont1=\eleveni\ifx\paradouze\oui\onzepts\scriptfont1=\ninei\elevenit\else
                        \douzepts\twelveit\fi}
 
   \def\fontetitredeuxb{\ifx\paradouze\oui\onzepts\eleventi\gponze\textfont1=\elevenib\scriptfont1=\nineib
                         \else\douzepts\twelveti\scriptfont1=\twelveib\scriptfont1=\tenib\gpdouze\fi}
                         
\def\fontetitredeuxl{\onzepts\textfont1=\elevenbf\scriptfont1=\ninebf\twelvebf}
  
\def\fontetitretrois{\textfont0=\elevenrm\scriptfont0=\eightrm\textfont1=\eleveni
                      \scriptfont1=\eighti\scriptscriptfont1=\sixi\elevenit}
                      
\def\fontetitrequatre{\textfont0=\elevenrm\scriptfont0=\eightrm\textfont1=\eleveni
                      \scriptfont1=\eighti\scriptscriptfont1=\sixi\elevenrm}
  
  \newcount\titreun\titreun=0
  \newcount\titredeux\titredeux=0
  \newcount\titretrois\titretrois=0
  \newcount\titrequatre\titrequatre=0
  \newcount\enonce\enonce=0
  
  \def\incr#1{\global\advance#1 by 1 {\the #1}}
  \def\avance#1{\global\advance#1 by 1}
  \def\init#1{\global#1=0}
  
  \long\def\Indentation#1#2{\setbox10=\hbox{\fontetitreun#1}
  	                    \ifdim\wd10 < 4mm
                         \setbox10=\hbox to 4mm{\box10\hfill}
                       \else\ifdim\wd10 < 6mm
                         \setbox10=\hbox to 6mm{\box10\hfill}
  	                    \else\ifdim\wd10 < 8mm
                         \setbox10=\hbox to 8mm{\box10\hfill}
                       \else\ifdim\wd10 < 12mm
                         \setbox10=\hbox to 12mm{\box10\hfill}
                       \fi\fi\fi\fi
                       \dimen10=\hsize
                       \advance \dimen10 by -\wd10
                       \noindent \box10 %
                       \ignorespaces
                       \hbox{\vtop{\hsize=\dimen10\raggedright\noindent\fontetitreun#2}}}

  \long\def\paraun#1{\removelastskip\par\medskip\goodbreak\vskip0pt plus.01\vsize\penalty-100
                \vskip0pt plus-.01\vsize
  	              \init{\titredeux}\ifnum\optionparag=1{\init\eqnumber\init\enonce}\else{}\fi
                  \goodbreak{\fontetitreun
  	                \Indentation{\incr{\titreun}.\ }{\fontetitreun #1\par}}\nobreak\medskip}

 %
 %
 \long\def\paraunc#1{\removelastskip\par\bigskip\goodbreak\vskip0pt plus.01\vsize\penalty-100
                \vskip0pt plus-.01\vsize
  	              \init{\titredeux}
                 \ifnum\optionparag=1{\init{\eqnumber}\init\enonce}\else{}\fi
                  \goodbreak
  	                {\parindent0mm\raggedcenter\fontetitreun\incr{\titreun}.\ 
                     \fontetitreun #1\par}\nobreak\medskip}
                     
\newtoks\titreunl
\titreunl={\ifnum\titreun=1{I}\fi%
\ifnum\titreun=2{II}\fi%
\ifnum\titreun=3{III}\fi%
\ifnum\titreun=4{IV}\fi%
\ifnum\titreun=5{V}\fi%
\ifnum\titreun=6{VI}\fi%
\ifnum\titreun=7{VII}\fi%
\ifnum\titreun=8{VIII}\fi%
\ifnum\titreun=9{IX}\fi%
\ifnum\titreun=10{X}\fi%
\ifnum\titreun=11{XI}\fi%
\ifnum\titreun=12{XII}\fi%
\ifnum\titreun=13{XIII}\fi%
}
\long\def\paraunl#1{\removelastskip\par\bigskip\bigskip\goodbreak\vskip0pt plus.01\vsize\penalty-100
                \vskip0pt plus-.01\vsize
  	              \init{\titredeux}\ifnum\optionparag=1{\init\eqnumber\init\enonce}\else{}\fi
                  \goodbreak{\fontetitreunl
  	                \Indentation{\global\advance\titreun by 1{\the\titreunl}.\ }{\fontetitreunl #1\par}}\nobreak\smallskip}

  
  \long\def\paradeux#1{\init{\titretrois}\vskip0pt plus.01\vsize\penalty-10
                \vskip0pt plus-.01\vsize\ifx \elie\oui\medskip\ifnum\titredeux>0\medskip\fi\fi
                 \Indentation{\fontetitredeux\the\titreun${\cdot}$\incr{\titredeux}.}
                              {\fontetitredeux\textfont1=\eleveni#1}\nobreak\par }
  
  \long\def\paradeuxb#1{\init{\titretrois}\vskip0pt plus.001\vsize\penalty-10
                \vskip0pt plus-.01\vsize{\ifx \elie\oui\medskip\ifnum\titredeux>0\medskip\fi\fi
                  \Indentation
  {\fontetitredeuxb\the\titreun${\cdot}$\incr{\titredeux}.}{ \fontetitredeuxb#1}}\nobreak
\smallskip}

\newtoks\titredeuxl
\titredeuxl={\ifnum\titredeux=1{A}\fi%
\ifnum\titredeux=2{B}\fi%
\ifnum\titredeux=3{C}\fi%
\ifnum\titredeux=4{D}\fi%
\ifnum\titredeux=5{E}\fi%
\ifnum\titredeux=6{F}\fi%
\ifnum\titredeux=7{G}\fi%
\ifnum\titredeux=8{H}\fi%
\ifnum\titredeux=9{I}\fi%
\ifnum\titredeux=10{J}\fi%
\ifnum\titredeux=11{K}\fi%
\ifnum\titredeux=12{L}\fi%
\ifnum\titredeux=13{M}\fi%
}
 \long\def\paradeuxl#1{\init{\titretrois}\vskip0pt plus.001\vsize\penalty-10
                \vskip0pt plus-.01
                \vsize \bigskip%
                  \Indentation
     {\fontetitredeuxl\global\advance\titredeux by 1
  \quad \the\titreunl${\cdot}$\the\titredeuxl.}{ \fontetitredeuxl#1}
  \removelastskip\nobreak\smallskip}
  

  \long\def\paratrois#1{\init{\titrequatre}\ifdim\lastskip<\smallskipamount
                \removelastskip\smallskip\fi
                 \vskip0pt plus.01\vsize\penalty-10
                  \vskip0pt
plus-.01\vsize{\ifx \elie\oui\ifnum\titretrois>0\medskip\fi\fi
\Indentation{\fontetitretrois\the\titreun${\cdot}$\the\titredeux${\cdot}$\incr{\titretrois}.\ }
  {\hskip0mm\baselineskip=14pt\fontetitretrois#1}\nobreak\smallskip}}
  
  
  \long\def\paratroisl#1{\init{\titrequatre}\ifdim\lastskip<\smallskipamount
                \removelastskip\fi
                 \vskip0pt plus.01\vsize\penalty-10
                  \vskip0pt
plus-.01\vsize\ifx \elie\oui\bigskip
\fi
\Indentation{\fontetitretrois\quad \quad \the\titreunl{${\cdot}$}\the\titredeuxl${\cdot}$\incr{\titretrois}.\ }
  {\hskip0mm\fontetitretrois#1}\nobreak\smallskip}


  \long\def\paraquatre#1{\ifdim\lastskip<\smallskipamount
                \removelastskip\smallskip\fi
                 \vskip0pt plus.01\vsize\penalty-10
                  \vskip0pt
                  plus-.01\vsize\par
 
\Indentation{\fontetitrequatre \the\titreun{${\cdot}$}\the\titredeux${\cdot}$\the\titretrois${\cdot}$\incr{\titrequatre}.\ }
{\hskip0mm\fontetitrequatre#1}\nobreak\smallskip}


\newtoks\titrequatrel
\titrequatrel={\ifnum\titrequatre=1{a}\fi%
\ifnum\titrequatre=2{b}\fi%
\ifnum\titrequatre=3{c}\fi%
\ifnum\titrequatre=4{d}\fi%
\ifnum\titrequatre=5{e}\fi%
\ifnum\titrequatre=6{f}\fi%
\ifnum\titrequatre=7{g}\fi%
\ifnum\titrequatre=8{h}\fi%
\ifnum\titrequatre=9{i}\fi%
\ifnum\titrequatre=10{j}\fi%
\ifnum\titrequatre=11{k}\fi%
\ifnum\titrequatre=12{l}\fi%
\ifnum\titrequatre=13{m}\fi%
}
\long\def\paraquatrel#1{\ifdim\lastskip<\smallskipamount
                \removelastskip\smallskip\fi
                 \vskip0pt plus.01\vsize\penalty-10
                  \vskip0pt
                  plus-.01\vsize{\bigskip
\Indentation{\global\advance\titrequatre by 1
\fontetitrequatre\quad \quad \quad \the\titreunl${\cdot}$\the\titredeuxl${\cdot}$\the\titretrois${\cdot}$\the\titrequatrel.\ }
{\hskip0mm\fontetitrequatre#1}\nobreak\smallskip}}

\ifx\optionkeys\oui
\def\drefun#1{\definexref{¤#1}{{\the\titreun}}{}} 
\def\drefdeux#1{\definexref{¤#1}{{\the\titreun}.{\the\titredeux}}{}}
\def\dreftrois#1{\definexref{¤#1}{{\the\titreun}.{\the\titredeux}.{\the\titretrois}}{}}
\else
\def\drefun#1{\definexref{prg#1}{{\the\titreun}}{}} 
\def\drefdeux#1{\definexref{prg#1}{{\the\titreun}.{\the\titredeux}}{}}
\def\dreftrois#1{\definexref{prg#1}{{\the\titreun}.{\the\titredeux}.{\the\titretrois}}{}}
\fi

%


  \long\def\propdeux#1#2#3#4{%
       \avance{\enonce}
       \leavevmode\edef\temp{#2}%
         \ifx\temp\empty 
          \else
           \definexref{#2}{#1~{\the\titreun.\the\enonce}}{enonces}
            \definexref{s#2}{{\the\titreun.\the\enonce}}{enonces}
             \fi
\smallskip
      \noindent{\bf#1\ {\bf\the\titreun.\the\enonce{#3}.}\enspace}{\sl#4\par}%
      \ifdim\lastskip<\medskipamount \removelastskip\penalty55\par \fi
   }

  \long\def\propun#1#2#3#4{%
      \avance{\enonce}
       \leavevmode\edef\temp{#2}%
        \ifx\temp\empty 
          \else
           \definexref{#2}{#1~{\the\enonce}}{enonces}
            \definexref{{s#2}}{{\the\enonce}}{enonces}
             \fi
   \par 
     \noindent{\bf#1\ {\bf\the\enonce{#3}.}\enspace}{\sl#4\par}%
     \ifdim\lastskip<\medskipamount \removelastskip\penalty55\medskip\fi
  }
  
  \long\def\prop#1#2#3#4{\ifnum\optionparag=1
                          \propdeux{#1}{#2}{\textfont1=\elevenib#3}{#4} \else\propun{#1}{#2}{\textfont1=\elevenib#3}{#4}\fi}

  \long\def\propt#1#2#3{\ifx\tpf\oui \prop{Th\'eo\-r\`eme}{#1}{#2}{#3}\par
                       \else\prop{Theorem}{#1}{#2}{#3}\par\fi}
  \long\def\Propt#1#2{\propt{#1}{}{#2}}
  \long\def\propl#1#2#3{\ifx\tpf\oui\prop{Lem\-me}{#1}{#2}{#3}\par
                         \else\prop{Lemma}{#1}{#2}{#3}\par\fi}
  
  \long\def\propc#1#2#3{\ifx\tpf\oui\prop{Corol\-laire}{#1}{#2}{#3}\par
                         \else\prop{Corollary}{#1}{#2}{#3}\par\fi}

  \long\def\propd#1#2#3{\ifx\tpf\oui\prop{D\'efi\-nition}{#1}{#2}{#3}\par
                       \else\prop{Definition}{#1}{#2}{#3}\par\fi} 
  
  \long\def\proptd#1#2#3{\ifx\tpf\oui\prop{Th\'eor\`eme et d\'efi\-nition}{#1}{#2}{#3}\par
                       \else\prop{Theorem and definition}{#1}{#2}{#3}\par\fi}


  
  \newcount\optionparag\optionparag=1
  
  \long\def\section#1#2{\ifnum\optionparag=1 \paraun{#2} 
                        \else\goodbreak{\fontetitreun
  	                \Indentation{#1.\ }{#2}}\nobreak\smallskip\fi}

  \def\eqconstruct#1{\ifnum\optionparag=1{\the\titreun\hbox{$\cdot$}#1}\else{#1}\fi}

  
  
  \def\numref{oui}  
  
  \newcount\mesref\mesref=0 
  \def\defbib#1{\ifx\numref\oui\global\advance\mesref by 1 \definexref{#1}{{\the
                 \mesref}}{}\else\definexref{#1}{#1}{}\fi}
  \def\bibtem#1{\defbib{#1}\item{\citer{#1}}}
  \def\citer#1{[\ref{#1}]}

  
  \font\seventeenmsa=msam10 at 17pt    
  \font\fourteenmsa=msam10 at 14pt
  \font\twelvemsa=msam10 at 12pt
  \font\tenmsa=msam10                 
  \font\ninemsa=msam10 at 9pt 
  \font\eightmsa=msam10 at 8pt 
  \font\sevenmsa=msam7 
  \font\sixmsa=msam10 at 6pt
  \font\fivemsa=msam5
  \newfam\msafam\textfont\msafam=\tenmsa\scriptfont\msafam=\sevenmsa\scriptscriptfont\msafam=\fivemsa
  
  \font\seventeenbb=msbm10 at 17pt     
  \font\fourteenbb=msbm10 at 14pt
  \font\twelvebb=msbm10 at 12pt
  \font\tenbb=msbm10                   
  \font\ninebb=msbm10 at 9pt 
  \font\eightbb=msbm10 at 8pt 
  \font\sevenbb=msbm7 
  \font\sixbb=msbm10 at 6pt
  \font\fivebb=msbm5 
  \newfam\bbfam\textfont\bbfam=\tenbb\scriptfont\bbfam=\sevenbb\scriptscriptfont\bbfam=\fivebb
  \def\bb{\fam\bbfam\tenbb}%

  \font\seventeenscaln=eusm10 at 17pt   
  \font\twelvescaln=eusm10 at 12pt
  \font\tenscaln=eusm10                
  \font\ninescaln=eusm10 scaled 900
  \font\eightscaln=eusm10 scaled 800
  \font\sevenscaln=eusm10 scaled 700
  \font\sixscaln=eusm10 scaled 600
   
  \newfam\scalnfam\textfont\scalnfam=\tenscaln\scriptfont\scalnfam=\sevenscaln\scriptscriptfont\scalnfam=\sixscaln
  \def\scaln{\fam\scalnfam\tenscaln}%
  \def\scal{\scaln}
  
  \font\tenscalb=eusb10                

  \font\sevenscalb=eusb10 scaled 700

  \newfam\scalbfam\textfont\scalbfam=\tenscalb\scriptfont\scalbfam=\sevenscalb
  %
  
  %
  %
  \font\fourteenrm=cmr12 scaled 1200
  \font\elevenrm=cmr10 at 11pt
  \font\twelverm=cmr12
  \font\ninerm=cmr9
  \font\eightrm=cmr8      
  \font\sevenrm=cmr7
  \font\sixrm=cmr6

  \font\seventeenpcap=cmcsc10 at 17pt
  \font\tenpcap=cmcsc10                        
  \font\ninepcap=cmcsc9
  \font\eightpcap=cmcsc8
  \font\sevenpcap=cmcsc10 scaled 700
  
  \newfam\pcapfam\textfont\pcapfam=\tenpcap\scriptfont\pcapfam=\sevenpcap
  \def\pcap{\fam\pcapfam\tenpcap}
  
  \font\seventeenrm=cmbx12 scaled 1400

  \font\fourteenbf=cmbx10 scaled 1400
  
  \font\twelvebf=cmbx12
  \font\elevenbf=cmbx10 at 11pt
  \font\ninebf=cmbx9  
  \font\eightbf=cmbx8
  \font\sixbf=cmbx6
  
  \font\tengot=eufm10                           
   
  \font\eightgot=eufm10 at 8truept 
  \font\sevengot=eufm7 
  \font\sixgot=eufm10 at 6 truept 
   
  \newfam\gotfam
  \textfont\gotfam=\tengot\scriptfont\gotfam=\sevengot\scriptscriptfont\gotfam=\sixgot
  %

  
  \def\tit{%
  \textfont0=\seventeenrm\scriptfont0=\tenrm\def\rm{\fam0\seventeenrm}%
  \textfont1=\seventeenib\scriptfont1=\twelveib%
  \textfont2=\seventeensy\scriptfont2=\twelvesy\scriptscriptfont2=\ninesy
  \textfont3=\seventeenex
  \textfont\itfam=\seventeenti
  \def\it{\fam\itfam\seventeenti}%
  \textfont\bbfam=\seventeenbb \scriptfont\bbfam=\twelvebb
  \def\bb{\fam\bbfam\seventeenbb}%
  \textfont\msafam=\seventeenmsa\scriptfont\msafam=\twelvemsa
  \textfont\scalnfam=\seventeenscaln
  \def\pcap{\fam\pcapfam\seventeenpcap}
  \normalbaselineskip=25pt\normalbaselines\rm}

  \font\seventeenti=cmbxti10 scaled 1680
  
  \font\fourteenti=cmbxti10 at 14pt
  
  \font\twelveti=cmbxti10 scaled 1200
  \font\eleventi=cmbxti10 at 11pt

  %
  %
  \font\twelveit=cmti12	
  \font\elevenit=cmti10 scaled 1100
  \font\nineit=cmti9
  \font\eightit=cmti8
  \font\sevenit=cmti7

  %
  %
 
 \font\seventeenib=cmmib10 scaled 1680
  \font\fourteenib=cmmib10 scaled 1400
  \font\twelveib=cmmib10 scaled 1200
  \font\elevenib=cmmib10 scaled 1100
  \font\tenib=cmmib10
\font\eightib=cmmib10 scaled 800
  \font\nineib=cmmib10 scaled 900
\font\sevenib=cmmib10 scaled 700
\font\sixib=cmmib10 scaled 600
\font\fiveib=cmmib10 scaled 500

\ifx\ITAN\oui
\else
\innernewfam\cmmibfam
\textfont\cmmibfam=\tenib
\scriptfont\cmmibfam=\sevenib
\scriptscriptfont\cmmibfam=\fiveib
\def\ib{\fam\cmmibfam\tenib}
\fi

  %
  %
  
  \font\eleveni=cmmi10 scaled 1100
  \font\ninei=cmmi9
  \font\eighti=cmmi8 
  \font\seveni=cmmi7 			                
  \font\sixi=cmmi6
  
  \font\ninesl=cmsl9                    
  \font\eightsl=cmsl8 
  \font\sevensl=cmsl10 at 7pt

  \font\ninett=cmtt9                    
  \font\eighttt=cmtt8
  \font\seventt=cmtt10 scaled 700

  \font\seventeensy=cmsy10 scaled 1680    
  \font\fourteensy=cmsy10 scaled 1400
  \font\twelvesy=cmsy10 scaled 1176
  
  \font\ninesy=cmsy9                      
  \font\eightsy=cmsy8
  \font\sixsy=cmsy6
  \font\seventeenex=cmex10 at 17pt
  \font\fourteenex=cmex10 at 14pt
  \font\twelveex=cmex10 at 12pt
  \font\nineex=cmex10 at 9pt
  \font\eightex=cmex10 at 8pt
  \font\sevenex=cmex10 at 7pt
  \font\sixex=cmex10 at 6pt
  \font\fiveex=cmex10 at 5pt
  
   
  \font\fourteengp=cmmi10 at 14pt
  
  \font\twelvegp=cmmib10 at 12pt
  \font\elevengp=cmmib10 at 11pt
  \font\tengp=cmmib10                          
  \font\ninegp=cmmib10 at 9pt 
  \font\eightgp=cmmib8 
   
  \font\sixgp=cmmib6


  \def\gponze{\textfont0=\elevenbf\scriptfont0=\eightbf\scriptscriptfont0=\sixbf
           \textfont1=\elevengp\scriptfont1=\eightgp\scriptscriptfont1=\sixgp}
  \def\gpdouze{\textfont0=\twelvebf\scriptfont0=\tenbf\scriptscriptfont0=\ninebf
           \textfont1=\twelvegp\scriptfont1=\tengp\scriptscriptfont1=\ninegp}        
  
 \def\gpquatorze{\textfont0=\fourteenbf\scriptfont0=\twelvebf\scriptscriptfont0=\elevenbf
           \textfont1=\fourteengp\scriptfont1=\twelvegp\scriptscriptfont1=\elevengp}

  
  \expandafter\chardef\csname pre amssym.def at\endcsname=\the\catcode`\@
  \catcode`\@=11
  \def\undefine#1{\let#1\undefined}
  \def\newsymbol#1#2#3#4#5{\let\next@\relax
   \ifnum#2=\@ne\let\next@\msafam@\else
   \ifnum#2=\tw@\let\next@\bbfam@\fi\fi
   \mathchardef#1="#3\next@#4#5}
  \def\mathhexbox@#1#2#3{\relax
   \ifmmode\mathpalette{}{\m@th\mathchar"#1#2#3}%
   \else\leavevmode\hbox{$\m@th\mathchar"#1#2#3$}\fi}
  \def\hexnumber@#1{\ifcase#1 0\or 1\or 2\or 3\or 4\or 5\or 6\or 7\or 8\or
   9\or A\or B\or C\or D\or E\or F\fi}
  
  \def\setboxz@h{\setbox\z@\hbox}
  \def\wdz@{\wd\z@}
  \def\boxz@{\box\z@}
  
  \edef\msafam@{\hexnumber@\msafam}
  \mathchardef\dabar@"0\msafam@39
  
  \edef\bbfam@{\hexnumber@\bbfam}
  \def\widehat#1{\setboxz@h{$\m@th#1$}%
   \ifdim\wdz@>\tw@ em\mathaccent"0\bbfam@5B{#1}%
   \else\mathaccent"0362{#1}\fi}
  \def\widetilde#1{\setboxz@h{$\m@th#1$}%
   \ifdim\wdz@>\tw@ em\mathaccent"0\bbfam@5D{#1}%
   \else\mathaccent"0365{#1}\fi}
  \newsymbol\leqq 1335          
  \newsymbol\leqslant 1336
  \newsymbol\lessgtr 1337       
  \newsymbol\backprime 1038     
  \newsymbol\risingdotseq 133A  
  \newsymbol\fallingdotseq 133B 
  \newsymbol\succcurlyeq 133C   
  \newsymbol\geqq 133D          
  \newsymbol\geqslant 133E
  \newsymbol\nmid 232D
  \newsymbol\nexists 2040
  \newsymbol\smallsetminus 2272
  \newsymbol\varnothing 203F
  
  \catcode`\@=\active

  \catcode`\@=11
  \newcount\typofr\typofr=1
  
  \catcode`\;=\active
  \def;{\ifnum\typofr=1\relax\ifhmode\ifdim\lastskip>\z@\unskip\fi
     \kern.2em\fi\string;\else\string;\fi}
  
  \catcode`\:=\active
  \def:{\ifnum\typofr=1\relax\ifhmode\ifdim\lastskip>\z@\unskip\fi
  \penalty\@M\ \fi\string:\else\string:\fi}
  
  \catcode`\!=\active
  \def!{\ifnum\typofr=1\relax\ifhmode\ifdim\lastskip>\z@\unskip\fi
     \kern.2em\fi\string!\else\string!\fi}
  
  \catcode`\?=\active
  \def?{\ifnum\typofr=1\relax\ifhmode\ifdim\lastskip>\z@\unskip\fi
     \kern.2em\fi\string?\else\string?\fi}

  \def\francais{\typofr=1\def\tpf{oui}}
  \def\anglais{\typofr=2\def\tpf{non}\def\english{oui}}
  \def\oui{oui}
  \francais
  
  \catcode`\@=12
  

\ifx\textures\oui
\def\raye #1|{\leavevmode\setbox1=\hbox{#1}%
\raise .5pt\hbox to \wd1{\xleaders\hbox{\rge{$ \sct / $}%
\kern 1pt}\hfill\kern -1pt }\kern -\wd1 \unhbox1\relax }

\def\barre #1|{\leavevmode\setbox1=\hbox{#1}%
\rlap{\Red\vrule height 2.4pt depth -1.2pt width \wd1}\Black \unhbox1\relax}
\else
\def\raye #1|{\leavevmode\setbox1=\hbox{#1}%
\raise .5pt\hbox to \wd1{\xleaders\hbox{\rge{$ \sct / $}%
\kern 1pt}\hfill\kern -1pt }\kern -\wd1 \unhbox1\relax }

\def\barre #1|{\leavevmode\setbox1=\hbox{#1}%
\rlap{\color{red}\vrule height 2.4pt depth -1.2pt width \wd1}\color{black} \unhbox1\relax}

\fi
  

  
  \def\og{\leavevmode\raise.24ex\hbox{$\scriptscriptstyle\langle\!\langle\>$}}    
  \def\fg{\leavevmode\raise.24ex\hbox{$\scriptscriptstyle\>\rangle\!\rangle$}}    
  \def\chv#1{\left\langle#1\right\rangle}

  \def\d{\,{\rm d}}
  \def\dd{{\rm d}}

  \def\z{{\bb Z}}
  \def\r{{\bb R}}
  
  \def\N{{\bb N}}

  \def\A{{\scal A}}
  \def\B{{\scal B}}

  \def\E{{\scal E}}
  \def\G{{\scal G}}

  \def\M{{\scal M}}
  
  \def\O{{\scal O}}
  \def\P{{\scaln P}}

  \def\T{{\scal T}}

  \def\frac#1#2{{#1\over #2}}
  \def\di#1#2{\sct#1\atop{\sct#2}}

  \def\numero{n$^{\rm o}\thinspace$}
\def\numeros{n$^{\rm os}\thinspace$}

  \def\qedbox{$\rlap{$\sqcap$}\sqcup$}           
  \def\qed{\nobreak\hfill\penalty250 \hbox{}\nobreak\hfill\qedbox\par }

  \def\numero{n$^{\rm o}\thinspace$}

  \def\¤{\S\thinspace}

  \def\¥{$\bullet$ }
  
  
  \def\e{{\rm e}}
  \def\mod{\mathop{\rm mod}\nolimits}
  \def\md#1#2{\equiv#1\,({\rm mod\,}#2)}
  \def\no#1{\Vert#1\Vert}

  \def\epsilon{\varepsilon}

  \def\phi{\varphi}
  \def\theta{\vartheta}
  \def\rho{\varrho}
  \def\dm{{\textstyle{1\over 2}}}
  \def\txt{\textstyle}
  \def\dsp{\displaystyle}
  \def\sct{\scriptstyle}
  \def\pf{\noi{\it Proof. }}
  \def\nid{\ifnum\typofr=1\par\noindent{\it D\'emonstration. }\else\pf\fi}
  \def\noi{\noindent}
  \def\rem{\ifnum\typofr=1\noi{\it Remarque.}\ \else\noi{\it Remark.}\ \fi}
  \def\rems{\ifnum\typofr=1\noi{\it Remarques.}\ \else\noi{\it Remarks.}\ \fi}

  \def\ov{\overline}
  \def\un{\underline}
  \def\sset{\smallsetminus}

  \def\pp{{\rm pp}}

  \def\1{{\bf 1}}
  \def\|{\Vert}

  \def\le{\leqslant}\def\leq{\leqslant}
  \def\ge{\geqslant}\def\geq{\geqslant}

  \def\ie{{i.e.\ }}
  \def\eg{{e.g.}}
  

  \def\log{\mathop{\rm log}\nolimits}
  \def\ft#1#2{{\txt{#1\over #2}}}



\def\Vbs#1{\bigg|#1\bigg|}


  \def\pmb#1{\setbox0=\hbox{#1}%
  \kern-.025em\copy0\kern-\wd0\kern.05em\copy0\kern-\wd0\kern-.025em\raise .0433em\box0 }

  
  \skewchar\eighti='177 \skewchar\sixi='177
  \skewchar\eightsy='60 \skewchar\sixsy='60
  
  \def\eightpoint{%
  \textfont0=\eightrm\scriptfont0=\sixrm\scriptscriptfont0=\fiverm
  \def\rm{\fam0\eightrm}%
  \textfont1=\eighti\scriptfont1=\sixi
  \scriptscriptfont1=\fivei\def\oldstyle{\fam1\seveni}%
  \textfont2=\eightsy\scriptfont2=\sixsy\scriptscriptfont2=\fivesy
  \textfont3=\eightex\scriptfont3=\sixex
  \textfont\itfam=\eightit
  \def\it{\fam\itfam\eightit}%
  \textfont\slfam=\eightsl
  \def\sl{\fam\slfam\eightsl}%
  \textfont\bbfam=\eightbb \scriptfont\bbfam=\sixbb\scriptscriptfont\bbfam=\fivebb
  \def\bb{\fam\bbfam\eightbb}%
  \textfont\msafam=\eightmsa\scriptfont\msafam=\sixmsa
  \textfont\scalnfam=\eightscaln
  \def\scaln{\fam\scalnfam\eightscaln}
  \textfont\ttfam=\eighttt
  \def\tt{\fam\ttfam\eighttt}%
\textfont\gotfam=\eightgot
  \textfont\bffam=\eightbf\scriptfont\bffam=\sixbf\scriptscriptfont\bffam=\fivebf
  \def\bf{\fam\bffam\eightbf}%
  \ifx\ITAN\oui\else\textfont\cmmibfam=\eightib
       \scriptfont\cmmibfam=\sixib
        \scriptscriptfont\cmmibfam=\fiveib
         \def\ib{\fam\cmmibfam\eightib}
   \fi
  \textfont\pcapfam=\eightpcap
  \def\pcap{\fam\pcapfam\eightpcap}
  \abovedisplayskip=2pt plus2pt minus 2pt
  \belowdisplayskip=2pt plus1pt minus 2pt
  \abovedisplayshortskip= 1pt plus 2pt minus 1pt
  \belowdisplayshortskip= 1pt plus 2pt minus 1pt
  \smallskipamount=2pt plus 1pt minus 2pt
  \medskipamount=3pt plus 2pt minus 2pt
  \bigskipamount=7pt plus 3pt minus 3pt
  \setbox\strutbox=\hbox{\vrule height 5pt depth 2pt width 0pt}%
  \normalbaselineskip=9pt\normalbaselines\rm}

  \def\({\left(}
  \def\){\right)}
  
  \def\footnoterule{\kern -2pt\hrule width 7truecm\kern 2.4pt}
  
  \def\xnotedef#1{\definexref{#1}{\noexpand\number\footnotenumber}{Note}}%

  
  
  \def\ninepoint{%
  \textfont0=\ninerm\scriptfont0=\sixrm\scriptscriptfont0=\fiverm
  \def\rm{\fam0\ninerm}%
  \textfont1=\ninei\scriptfont1=\sixi
  \scriptscriptfont1=\fivei\def\oldstyle{\fam1\ninei}%
  \textfont2=\ninesy\scriptfont2=\sixsy\scriptscriptfont2=\fivesy
  \textfont3=\nineex\scriptfont3=\sixex
  \textfont\itfam=\nineit
  \def\it{\fam\itfam\nineit}%
  \textfont\slfam=\ninesl
  \def\sl{\fam\slfam\ninesl}%
  \textfont\bbfam=\ninebb\scriptfont\bbfam=\sixbb\scriptscriptfont\bbfam=\fivebb
  \def\bb{\fam\bbfam\ninebb}%
  \textfont\msafam=\ninemsa\scriptfont\msafam=\sixmsa\scriptscriptfont\msafam=\fivemsa
  \textfont\scalnfam=\ninescaln
  \def\scaln{\fam\scalnfam\ninescaln}
  \textfont\ttfam=\ninett
  \def\tt{\fam\ttfam\ninett}%
  \textfont\bffam=\ninebf\scriptfont\bffam=\sixbf\scriptscriptfont\bffam=\fivebf
  \def\bf{\fam\bffam\ninebf}%
  \abovedisplayskip=3pt plus2pt minus 2pt
  \belowdisplayskip=3pt plus1pt minus 2pt
  \abovedisplayshortskip= 2pt plus 2pt minus 1pt
  \belowdisplayshortskip= 2pt plus 2pt minus 1pt
  \smallskipamount=2pt plus 1pt minus 2pt
  \medskipamount=3pt plus 2pt minus 2pt
  \bigskipamount=7pt plus 3pt minus 3pt
  \setbox\strutbox=\hbox{\vrule height 5pt depth 2pt width 0pt}%
  \normalbaselineskip=10.5pt plus.3pt minus.3pt\normalbaselines\rm}

  \def\sevenpoint{%
  \textfont0=\sevenrm\scriptfont0=\sixrm\scriptscriptfont0=\fiverm
  \def\rm{\fam0\sevenrm}%
  \textfont1=\seveni\scriptfont1=\sixi
  \scriptscriptfont1=\fivei\def\oldstyle{\fam1\seveni}%
  \textfont2=\sevensy\scriptfont2=\sixsy\scriptscriptfont2=\fivesy
  \textfont3=\sevenex\scriptfont3=\fiveex
  \textfont\itfam=\sevenit
  \def\it{\fam\itfam\sevenit}%
  \textfont\slfam=\sevensl
  \def\sl{\fam\slfam\sevensl}%
  \textfont\bbfam=\sevenbb \scriptfont\bbfam=\sixbb\scriptscriptfont\bbfam=\fivebb
  \def\bb{\fam\bbfam\sevenbb}%
  \textfont\msafam=\sevenmsa\scriptfont\msafam=\sixmsa
  \textfont\scalnfam=\sevenscaln
  \def\scaln{\fam\scalnfam\sevenscaln}
  \textfont\bffam=\sevenbf\scriptfont\bffam=\sixbf\scriptscriptfont\bffam=\fivebf
  \def\bf{\fam\bffam\sevenbf}%
  \textfont\ttfam=\seventt
  \abovedisplayskip=2pt plus2pt minus 2pt
  \belowdisplayskip=2pt plus1pt minus 2pt
  \abovedisplayshortskip= 1pt plus 2pt minus 1pt
  \belowdisplayshortskip= 1pt plus 2pt minus 1pt
  \smallskipamount=2pt plus 1pt minus 2pt
  \medskipamount=3pt plus 2pt minus 2pt
  \bigskipamount=7pt plus 3pt minus 3pt
  \setbox\strutbox=\hbox{\vrule height 5pt depth 2pt width 0pt}%
  \normalbaselineskip=9pt\normalbaselines\rm}

 \def\onzepts{%
 \textfont0=\elevenrm\scriptfont0=\ninerm
 \textfont1=\elevenib\scriptfont1=\ninei}

\def\douzepts{%
  \textfont0=\twelverm\scriptfont0=\tenrm\def\rm{\fam0\twelverm}%
  \textfont1=\twelveib\scriptfont1=\teni%
  \textfont2=\twelvesy\scriptfont2=\tensy\scriptscriptfont2=\eightsy
  \textfont3=\twelveex
  \textfont\itfam=\twelveti
  \def\it{\fam\itfam\twelveti}%
  \textfont\bffam=\twelvebf\scriptfont\bffam=\tenbf\scriptscriptfont\bffam=\eightbf
  \def\bf{\fam\bffam\twelvebf}%
  \textfont\bbfam=\twelvebb \scriptfont\bbfam=\tenbb
  \def\bb{\fam\bbfam\twelvebb}%
  \textfont\msafam=\twelvemsa\scriptfont\msafam=\tenmsa
  \textfont\scalnfam=\twelvescaln
  \normalbaselineskip=15pt\normalbaselines\rm}

\def\quatorzepts{%
  \textfont0=\fourteenrm\scriptfont0=\twelverm\def\rm{\fam0\fourteenrm}%
  \textfont1=\fourteenib\scriptfont1=\twelveib%
  \textfont2=\fourteensy\scriptfont2=\twelvesy\scriptscriptfont2=\tensy
  \textfont3=\fourteenex
  \textfont\itfam=\fourteenti
  \def\it{\fam\itfam\fourteenti}%
  \textfont\bffam=\fourteenbf\scriptfont\bffam=\twelvebf\scriptscriptfont\bffam=\tenbf
  \def\bf{\fam\bffam\fourteenbf}%
  \textfont\bbfam=\fourteenbb \scriptfont\bbfam=\twelvebb
  \def\bb{\fam\bbfam\fourteenbb}%
  \textfont\msafam=\fourteenmsa\scriptfont\msafam=\twelvemsa
  \textfont\scalnfam=\twelvescaln
  \normalbaselineskip=18pt\normalbaselines\rm}


\def\AA{{\it Acta Arith.}}

\def\JLMS{{\it J. London Math. Soc.}}

\def\PLMS{{\it Proc. London Math. Soc.}}

\def\picture #1 by #2 (#3){\leavevmode\vbox to #2{
     \hrule width #1 height 0pt depth 0pt
      \vfill
       \special{picture #3}}}

\def\illustration #1 by #2 (#3) scaled #4{\dimen1=#2
  \divide\dimen1 by 1000
  \multiply\dimen1 by #4
  \vtop to \dimen1{\dimen1=#1
  \divide\dimen1 by 1000
  \multiply\dimen1 by #4
  \hsize=\dimen1\vss
  \noindent\special{illustration #3 scaled #4}}}

\ifx\optionkeymacros\undefined\else \fi

\catcode`\Œ=\active\defŒ{{\aa}}       
\catcode`\º=\active\defº{\int}        
\catcode`\=\active\def{\c c}        
\catcode`\¶=\active\def¶{\partial}    
\catcode`\Ä=\active\defÄ{\oint}       
\catcode`\Æ=\active\defÆ{\triangle}   
\catcode`\Â=\active\defÂ{\neg}        
\catcode`\µ=\active\defµ{\mu}         
\catcode`\¿=\active\def¿{{\o}}        
\catcode`\¹=\active\def¹{\pi}         
\catcode`\Ï=\active\defÏ{{\oe}}       
\catcode`\§=\active\def§{{\ss}}       
\catcode`\ =\active\def {\dagger}     
\catcode`\Ã=\active\defÃ{\sqrt}       
\catcode`\·=\active\def·{\Sigma}      
\catcode`\Å=\active\defÅ{\approx}     
\catcode`\½=\active\def½{\Omega}      
\catcode`\£=\active\def£{{\it\$}}     
\catcode`\°=\active\def°{\infty}      
\catcode`\¤=\active\def¤{{\S}}        
\catcode`\¦=\active\def¦{{\P}}        
\catcode`\¥=\active\def¥{\bullet}     
\catcode`\»=\active\def»{\leavevmode\raise.585ex\hbox{\b a}}      
\catcode`\¼=\active\def¼{\leavevmode\raise.6ex\hbox{\b o}}        
\catcode`\­=\active\def­{\not=}       
\catcode`\²=\active\def²{\leq}        
\catcode`\³=\active\def³{\geq}        
\catcode`\Ö=\active\defÖ{\div}        
\catcode`\É=\active\defÉ{{\dots}}     
\catcode`\¾=\active\def¾{{\ae}}       
\catcode`\Ç=\active\defÇ{\og}         
\catcode`\Ò=\active\defÒ{``}          
\catcode`\Á=\active\defÁ{!`}          
\catcode`\¢=\active\def¢{\rlap/c}     
\catcode`\Ô=\active\defÔ{`}           
\catcode`\Õ=\active\defÕ{'}           


\catcode`\=\active\def{{\AA}}       
\catcode`\'=\active\def'{\c C}        
\catcode`\¯=\active\def¯{{\O}}        
\catcode`\¸=\active\def¸{\Pi}         
\catcode`\Î=\active\defÎ{{\OE}}       
\catcode`\®=\active\def®{{\AE}}       
\catcode`\×=\active\def×{\diamond}    
\catcode`\¡=\active\def¡{\accent'27}  
\catcode`\Ó=\active\defÓ{''}          
\catcode`\±=\active\def±{\pm}         
\catcode`\È=\active\defÈ{\fg}         
\catcode`\À=\active\defÀ{?`}          
\catcode`\Ð=\active\defÐ{--}          
\catcode`\Ñ=\active\defÑ{---}         


\catcode`\Š=\active\defŠ{\"a}        
\catcode`\'=\active\def'{\"e}        
\catcode`\•=\active\def•{\"{\i}}     
\catcode`\š=\active\defš{\"o}        
\catcode`\Ÿ=\active\defŸ{\"u}        
\catcode`\Ø=\active\defØ{\"y}        
\catcode`\å=\active\defå{\^A}        
\catcode`\€=\active\def€{\"A}        
\catcode`\…=\active\def…{\"O}        
\catcode`\†=\active\def†{\"U}        
\catcode`\‡=\active\def‡{\'a}        
\catcode`\Ž=\active\defŽ{\'e}        
\catcode`\'=\active\def'{\'{\i}}     
\catcode`\—=\active\def—{\'o}        
\catcode`\œ=\active\defœ{\'u}        
\catcode`\ƒ=\active\defƒ{\'E}        
\catcode`\æ=\active\defæ{\^E}        
\catcode`\é=\active\defé{\`E}        %
\catcode`\ˆ=\active\defˆ{\`a}        
\catcode`\=\active\def{\`e}        
\catcode`\"=\active\def"{\`{\i}}     
\catcode`\˜=\active\def˜{\`o}        
\catcode`\=\active\def{\`u}        
\catcode`\Ë=\active\defË{\`A}        
\catcode`\‹=\active\def‹{\~a}        
\catcode`\–=\active\def–{\~n}        
\catcode`\›=\active\def›{\~o}        
\catcode`\Ì=\active\defÌ{\~A}        
\catcode`\"=\active\def"{\~N}        
\catcode`\Í=\active\defÍ{\~O}        
\catcode`\‰=\active\def‰{\^a}        
\catcode`\=\active\def{\^e}        
\catcode`\"=\active\def"{\^{\i}}     
\catcode`\™=\active\def™{\^o}        
\catcode`\ž=\active\defž{\^u}        

\let\optionkeymacros\null

\dimstand
\optionparag=2
\anglais
\def\ITAN{oui}

    \font\tenrsfs=rsfs7 at 10pt

    \font\sevenrsfs=rsfs7
    \font\sixrsfs=rsfs7 at 6pt

\newfam\rsfsfam\textfont\rsfsfam=\tenrsfs\scriptfont\rsfsfam=\sevenrsfs\scriptscriptfont\rsfsfam=\sixrsfs
    \def\rsfs{\fam\rsfsfam\tenrsfs}%
\def\G{{\rsfs G}}
\def\paradouze{oui}
\def\auteur{GŽrald Tenenbaum}
\def\titrart{Some of Erd\H os' unconventional problems in number theory, thirty-four years later\anote{*}{We include here some corrections with respect to the published version.\hfill}}
\hautspages{\auteur}{\titrart}
{\leftskip-10mm\obeylines
 L. Lov‡sz  I.Z. Ruzsa,  V. T. S—s (eds), 
{\it Erd\H os Centennial volume}, 
Bolyai Society Mathematical Studies 
25 (2013), 651-681.\par }
\titrecentre{\titrart}
\bigskip\medskip
\centerline {\auteur} 
\bigskip\bigskip\bigskip
There are many ways to recall Paul Erd\H os' memory and his special way of doing mathematics. Ernst Straus described him as ``the prince of problem solvers and the absolute monarch of problem posers". Indeed, those mathematicians who are old enough to have attended some of his lectures will remember that, after his talks, chairmen used to slightly depart from standard conduct, not asking if there were any questions but if there were any answers.
\par 
In the address that he forwarded to Mikl—s Simonovits for Erd\H os' funeral, Claude Berge mentions a conversation he had with Paul in the gardens of the Luminy Campus, near Marseilles, in September 1995. After Paul's opening lecture for this symposium on Combinatorics, Berge asked him to specify his beauty criteria for a conjecture in discrete mathematics. Erd\H os mainly retained the following five:
\par 
{\leftskip 5mm \rightskip 5mm
(i) The {\it simplicity} of the statement;\par 
(ii) The expected {\it difficulty} of the  solution (which Paul liked to measure in dollars);
\par 
(iii) The {\it posterity} of the subsequent theorem, \ie the set of results arising either directly from the solution of from the methods designed to obtain it;\par 
(iv) The {\it future} of the path opened by the problem, which I would rather call the set of {\it descendants} of  the problem, in other words the family of new questions opened up by the statement or the solution of the conjecture;\par 
(v) The {\it intuitive representability} of the specific mathematical property that is being dealt with.
\par}
Apart, perhaps, the last, for which an adequate transposition should be described with further precision, these criteria are equally relevant to a classification for a conjecture in analytic and/or elementary number theory.
\par 
My purpose here mainly consists in illustrating these criteria by revisiting some of the problems stated by Erd\H os in his profound article \citer{Er79}.
\par 
Aside from updating the status of a number of interesting questions, my hope is to convince the reader that Erd\H os' conjectures, although stated in a condensed and seemingly particular form, were problematics rather than problems. Day after day, year after year, each of his questions appears, in the light of discussions and partial progress, as a node in a gigantic net, designed not for a single prey but for a whole species. 
\par 
In the sequel of this paper, quotes from the article \citer{Er79} are set in italics. I took liberties to correct obvious typographic errors and to slightly modify some notations in order to fit with subsequent works.
Erd\H os' paper starts with the following.
\par \smallskip
{\it First of all I state a very old conjecture of mine: the density of integers $n$ which have two divisors $d_1$ and $d_2$ satisfying $d_1<d_2<2d_1$ is 1. I proved long ago \citer{Er48} that the density of these numbers exists but I have never been able to prove that it is 1. I claimed \citer{Er64} that I proved that almost all integers $n$ have two divisors
$$d_1<d_2<d_1\big\{1+(\e/3)^{(1-\eta)\log\log n}\big\}\eqdef{conjdd'}$$
and that \eqref{conjdd'} is best possible, namely it fails if $1-\eta$ is replaced by $1+\eta$. R.R. Hall and I confirmed this later statement but unfortunately we cannot prove \eqref{conjdd'}. We are fairly sure that \eqref{conjdd'} is true and perhaps it is not hopeless to prove it by methods of probabilistic number theory that are at our disposal.}
\smallskip\goodbreak
This is an edifying example of a conjecture meeting the  above five requirements. However, before  elaborating on this, it may be worthwhile try understanding the process that led Erd\H os to this simple and deep statement.
\par \smallskip
An integer $n$ is called {\it perfect} if it is equal to the sum of its proper divisors. Thus $6=1+2+3$ and $28=1+2+4+7+14$ are perfect. In modern notation, a perfect integer $n$ satisfies $\sigma(n)=2n$ where $\sigma(n)$ stands for the sum of all divisors. This is an interesting formulation since $\sigma(n)$ is a multiplicative function of $n$.   In the third century before our era, Euclid proved (IX.36) that $2^{p-1}(2^p-1)$ is perfect whenever $2^p-1$ is prime, which of course implies that $p$ itself is prime.\par 
An integer $n$ is called {\it abundant} if $\sigma(n)>2n$.
In the early thirties, in a book on number theory, Erich Bessel-Hagen asks whether abundant integers have a natural density. Davenport \citer{Da33}, Chowla \citer{Ch34}, Erd\H os \citer{Er34} and Behrend \citer{Be35} all gave, independently, a positive answer. All proofs, except that of Erd\H os, rest on the method of (real or complex) moments. Erd\H os attacks the problem from another viewpoint: primitive abundant numbers, \ie abundant numbers having no abundant proper divisor. Writing $f(n)$ for $\sigma(n)/n$, any primitive abundant integer $n$ satisfies
$$2\leqslant f(n)\leqslant f(n/p)f(p)<2(1+1/p)$$
whenever $p|n$. Since the largest prime factor of $n$ is usually large, this restricts the cardinality of primitive abundant numbers not exceeding $x$, which can be shown to be $o\big(x/(\log x)^2\big)$. The proof is then completed by noticing that, if we write $$\M(\A):=\{ma:a\in\A,\,m\geqslant 1\}$$ for the so-called {\it set of multiples} of the set $\A$ and $\dd$, $\ov\dd$, $\underline\dd$ for natural, upper and lower density respectively, then
$$\dd\M(\A_T)\leqslant \underline\dd\M(\A)\leqslant \ov\dd\M(\A)\leqslant \dd\M(\A_T)+\sum_{\di{a>T}{a\in\A}}{1\over a}$$
holds for any integer sequence $\A$ such that $\sum_{a\in\A}1/a<\infty$, with $\A_T:=\A\cap[1,T]$.
\par 
This was the starting point of the fruitful concept of set of multiples.\par 
It was once suspected that any set of multiples should have a natural density. However, Besicovitch \citer{Bes34} soon disproved this conjecture by showing that
$$\liminf_{T\to\infty}\d\M(]T,2T])=0.\eqdef{Bes}$$
Indeed, it is easy to deduce from this that, given any $\varepsilon>0$ and a sequence $\{T_j\}_{j=0}^{\infty}$ increasing sufficiently fast, then $\A:=\cup_j]T_j,2T_j]$ satisfies $\underline\dd\M(\A)<\varepsilon$, $\ov\dd\M(\A)\geqslant \dm$.\par 
The reader might ask at this stage: interesting indeed, but how does this link to \eqref{conjdd'}? We still need a few more steps inside Erd\H os' peculiar way of thinking. \par 
It is one of the marks of the great: not to accept an obstruction before understanding it completely. This holds outside of mathematics as well as inside. Erd\H os did not accept Besicovitch's counter-example for itself and continued the quest.\par  First \citer{Er36}, he improved  \eqref{Bes} to the optimal
$$\lim_{T\to\infty}\d\M(]T,T^{1+\varepsilon_T}])=0\eqdef{E/B}$$
provided $\varepsilon_T\to0$ as $T\to\infty$. \par 
With this new, crucial piece of information, he progressed in two connected directions: first, to show, with Davenport \citer{DE37} --- see also \citer{DE51} for another, very interesting proof --- that any set of multiples has a logarithmic density, equal to its lower asymptotic density,\note{We shall make use of this extra information later on.} and, second,  to show \citer{Er48}\note{See \citer{EHT94} for a short proof.} that Besicovitch-type constructions are essentially the only obstacles to the existence of $\dd\M(\A)$: writing $d_1(n,\A):=\inf\{d|n:d\in\A\}$ with the convention that $d_1(n,\A)=\infty$ whenever $n\not\in\M(\A)$, a necessary and sufficient condition that $\M(\A)$ has a natural density is
$$\lim_{\varepsilon\to0}\ov \dd\{n\geqslant 1:n^{1-\varepsilon}<d_1(n,\A)\leqslant n\}=0.\eqdef{CNSdMA}$$ 
\par 
Now, consider the set $$\E:=\{m\in\N^*: m=dd', d<d'<2d\}.\eqdef{dd'}$$ Then $n^{1-\varepsilon}<d_1(n,\E)\leqslant n$ plainly implies that $n$ has a divisor in $]n^{1/2-\varepsilon},n^{1/2}]$ and it is easy to deduce \eqref{CNSdMA} from \eqref{E/B}.\par 
So we now know that the set of integers with two close divisors has a natural density. (By `close' we mean here that the ratio of the two divisors should lie in $]1,2[$.) Moreover, as seen above, the existence property follows in a natural way from the theory of sets of multiples: the sequence~$\E$ defined above is one of simplest examples one can think of that meets the criterion \eqref{CNSdMA}.
\par 
But what should the density be? Erd\H os stated, as early as 1948 (and probably much before) \citer{Er48}, that this density should be equal to 1. Here again, a seemingly anecdotal conjecture is actually based on a profound assumption---any answer to it, positive or negative, is bound to enlighten our understanding of the multiplicative structure of integers. 
\par 
 Let us make the convention to use the suffix $\pp$ to indicate that a relation holds on a set of asymptotic density 1. As we shall see later in this paper, Erd\H os had known for long that sufficiently far prime factors behave almost independently \pp. Specifically, if we denote~by
$$\{p_j(n)\}_{j=1}^{\omega(n)}\eqdef{sfp}$$
the increasing sequence of distinct prime factors of an integer $n$ and if we write $$U_j(n):=\{\log_2p_j(n)-j\}/\sqrt{j},\eqdef{Ujn}$$ then, to a first approximation, $U_j(n)$ and $U_h(n)$ resemble independent Gaussian random variables \pp\ provided that $j/h\to\infty$. (Here and in the sequel, we let $\log _k$ denote the $k$-fold iterated logarithm.) Having this in mind, it is reasonable to believe that, in first approximation, the quantities $\log (d'/d)$ are evenly distributed \pp\ in the interval $[-\log n,\log n]$. Since these quantities are $3^{\omega(n)}$ in number, we deduce from the Hardy--Ramanujan estimate $\omega(n)\sim\log_2n$ \pp\ that the smallest of these numbers should be of size $(\log n)^{1-\log 3+o(1)}$ \pp. \par 
This is, perhaps no more, certainly no less, what is hidden behind conjecture \eqref{conjdd'}.
\par \medskip
This conjecture, which is now a theorem, due to Erd\H os--Hall \citer{EH79} for the lower bound and to Maier--Tenenbaum \citer{MT84} for the upper bound, has had a wide posterity and many descendants. 
\par 
In his doctoral dissertation supervised by the author \citer{St92}, Stef proves that the number $R_x$ of exceptional integers not exceeding $x$ and which do not belong to $\M(\E)$ satisfies
$$x/(\log x)^{\beta+o(1)}\ll R_x\ll x\e^{-c\sqrt{\log_2 x}}\eqdef{encStef}$$
for a suitable constant $c>0$, with $\beta=1-(1+\log_23)/\log3\approx0,00415.$ These are the best known estimates to date.
\par \goodbreak
To the chapter of posterity certainly belong all results involving the still mysterious Erd\H os--Hooley Delta-function and the so-called propinquity functions
$$E_r(n):=\min_{1\leqslant j\leqslant \tau(n)-r}\log\{d_{j+r}(n)/d_j(n)\}\qquad (r\geqslant 1),$$
where $\{d_j(n)\}_{j=1}^{\tau(n)}$ stands for the increasing sequence of the divisors of an integer $n$.\par \goodbreak
One of the most recent achievements in this direction is a very precise confirmation of the heuristic principle leading to \eqref{conjdd'}, as described above: Raouj, Stef and myself prove in \citer{RST11} that 
$$E_1(n)={\log n\over 3^{\omega(n)}}(\log_2n)^{\vartheta_n}\qquad \pp,$$
where $-5\leqslant \vartheta_n\leqslant 10$. Many more precise and connected results are actually proved in~\citer{RST11}.\par \goodbreak
The situation is much less satisfactory regarding the functions $E_r$ when $r\geqslant 2$, for which the precise \pp\ behaviour is still unknown. Using techniques similar to that of the proof of
theorem 3 of
\citer{ET89}, it can be shown that 
$$E_2(n)>(\log n)^{-\gamma_2+o(1)} \qquad \pp$$
for some $\gamma_2<\log 3-1$. Moreover, the methods and results of \citer{MT11} yield
$$E_r(n)\leqslant (\log n)^{-\beta_r+o(1)}\qquad \pp,$$
with $$\beta_r:={(\log 3-1)^m\over (\log 3-1/3)^{m-1}},\quad2^{m-1}<r+1\leqslant 2^m.$$
Thus, we have
$$\beta_1=\log 3-1\approx 0.09861,\quad \beta_2=\beta_3\approx0.01271,\quad \beta_r\approx0.00164\quad(4\leqslant r\leqslant 7). $$
Also, it is proved in \citer{MT11} (th. 1.1) that $E_r(n)>\tau(n)^{-1/r+o(1)}$ holds \pp\ uniformly in $r\geqslant 1$, and thus
$$E_r(n)=1/(\log n)^{o(1)}\qquad \pp\quad(r=r(n)\to\infty),$$
a result which might look surprising at first sight.\par 
We conjecture the existence of a strictly decreasing sequence $\{\alpha_r\}_{r=1}^{\infty}$ such that
$$E_r(n)= (\log n)^{-\alpha_r+o(1)}\qquad \pp.$$
It is particularly irritating, for instance, to be unable to find a better \pp\ upper bound for $E_2(n)$ than for $E_3(n)$.
\par 
We also mention as a posterity result the proof by Raouj \citer{Ra95} of Erd\H os' conjecture asserting that $$\dd\M\big(\cup_{d|n}]d,2d]\big)=1+o(1) \qquad \pp.$$ This is established
 in the following fairly strong (and optimal) form. Put $\lambda^*:=\log 4-1$ and $\delta_n:=\dd\M\big(\cup_{d|n}]d,(1+1/(\log n)^\lambda)d]\big)$. Then
$$\eqalign{{1\over (\log n)^{F(\lambda)+o(1)}}<1-&\delta_n<\e^{-c_\lambda\sqrt{\log n}}\qquad\qquad(0\leqslant \lambda<\lambda^*) \cr
\delta_n&=(\log n)^{-F(\lambda)+o(1)}\qquad (\lambda>\lambda^*)}\qquad \pp,$$ 
where $F(\lambda):=\beta\log \beta-\beta+1$ with $\beta:=-1+(1+\lambda)/\log 2$ if $\lambda\leqslant 3\log 2-1$, and $F(\lambda):=\lambda-\log 2$ if $\lambda>3\log 2-1$.
\medskip
The Erd\H os--Hooley function is defined as
$$\Delta(n):=\sup_{u\in\r}\sum_{\di{d|n}{\e^u<d\leqslant \e^{u+1}}}1\qquad (n\geqslant 1).$$
It first appears  (implicitly) in \citer{Er74} and (explicitly) in \citer{EN75}, \citer{EN76} in the early seventies. It was next studied by Hooley \citer{Ho79} with the aim of developing a variety of applications to several branches of number theory. 
\par \goodbreak
The ratio $\Delta(n)/\tau(n)$ has an immediate probabilistic interpretation: with LŽvy's 1937 definition, it is the value at 1 of the concentration function of the random variable $D_n$ taking the values $\log d$ $(d|n)$ with uniform probability $1/\tau(n)$. It is noteworthy to state here that $D_n=\sum_{p^\nu\|n}D_{p^\nu}$ where the $D_{p^\nu}$ are independent. \par \smallskip
If we replace the factor 2
 by $\e$, which is irrelevant to all intents and purposes, Erd\H os' initial conjecture $$\dd\M(\E)=1\eqdef{conjE}$$ is equivalent to the statement that $$\Delta(n)>1\qquad \pp,\eqdef{Dn1}$$
so that \eqref{encStef} provides quantitative estimates for the number of exceptions.\par \goodbreak
The best $\pp$-bounds to date for the $\Delta$-function appear in a joint article with  Maier \citer{MT11}. We prove that
$$(\log_2n)^{\gamma+o(1)}<\Delta(n)<(\log_2n)^{\log 2+o(1)}\qquad \pp,$$
 where the exponent $\gamma:=(\log 2)/\log\big({1-1/\log 27\over 1-1/\log 3}\big)\approx0.33827$
is conjectured to be optimal.
\par 
To show the existence and determine the value of the exact exponent is a challenging problem in probabilistic number theory. There is no doubt that such a result would imply deeper ideas on the structure of the set of divisors of a normal integer.
\par 
However, as shown by Hooley in \citer{Ho79}, it is mainly information on the average order
$$s(x):={1\over x}\sum_{n\leqslant x}\Delta(n) $$
that has applications to other arithmetical topics such as Waring-type problems \citer{RCV85}, Diophantine approximation \citer{Ho79}, \citer{Te86}, and Chebyshev's problem on the greatest prime factor of polynomial sequences \citer{Te90}. It is thus proved in \citer{Te90}, as a consequence of an average estimate for a variant of $s(x)$, that, for any $\alpha<2-\log 4\approx0.61370$, the bound 
$$P^+\Big(\prod_{n\leqslant x}F(n)\Big)>x\,\e^{(\log x)^\alpha}\qquad (x>x_0(F))$$
holds for any irreducible polynomial $F(X)\in\z[X]$ with degree $>1$. This is currently the best available result valid for polynomials of arbitrary degree.
Here and in the sequel $P^+(m)$ denotes the largest prime factor of the integer $m$ with the convention that $P^+(1)=1$. \par 
 Established in \citer{HT82} and \citer{Te85}, the
best bounds for $s(x)$ at the time of writing~are
$$\log_2x\ll s(x)\ll \e^{c\sqrt{\log_2x\log_3x}}\qquad (x\to\infty) \eqdef{encormoy}$$
where $c$ is a suitable constant.
See \citer{HT88},
\citer{Te86} and, for instance, \citer{Ro11} for further references and descriptions on this question.
\par 
Still in the area of descendants of the conjecture \eqref{conjdd'}, we mention the recent paper \citer{BT12} in collaboration with La Bret{\`e}che and where sharp, weighted average bounds are given for functions of the type $$\Delta(n,f):=\sup_{u\in\r,\,0\leqslant
v\leqslant
1}\Vbs{\sum_{\di{\scriptscriptstyle d|n}{\scriptscriptstyle\e^{u}<d\leqslant
\e^{u+v}}}f(d)}\eqdef{Dnf}$$
where $f$ is an  oscillating function, typical cases being those of a non principal Dirichlet character or of the Mšbius function. All suitably weighted finite integral, even moments are also studied.  This is the key step to the proof, given in \citer{BT12+},
of Manin's conjecture, in the strong form conjectured by Peyre and with effective remainder term, for all Ch‰telet surfaces.
\par \goodbreak
Maier established in \citer{Ma87} normal upper and lower bounds for \eqref{Dnf} in the case $f=\mu$, the Mšbius function, and his method is equally applicable in the case $f=\chi$, a real, non principal Dirichlet character.\par 
Short averages have also been investigated, by Nair--Tenenbaum \citer{NT98}, Henriot \citer{He12}, and La Bretche--Tenenbaum \citer{BT12a}. These may have numerous, sometimes surprising applications. For instance, writing $\chv t$ for the fractional part of a real number $t$, we have \citer{NT98}, for any given $\varepsilon>0$,
$$\sup_{D\geqslant 1}\Vbs{\sum_{D\leqslant d\leqslant 2D}\chv{x+y\over d}-\chv{x\over d}}\ll
y(\log x)^{o(1)}\quad (x^\varepsilon\leqslant y\leqslant x), $$
a bound which known exponential sums methods, by far, will fail to meet.\par 
This ends our comments and update on conjecture \eqref{conjdd'}. 
\goodbreak
\medskip
The next problem in \citer{Er79} is described as follows.\par 
{\it Denote by $\tau^+(n)$ the number of integers $k$ for which $n$ has a divisor $d$ satisfying $2^k<d\leqslant 2^{k+1}$. I conjecture that for almost all $n$
$$\tau^+(n)/\tau(n)\to0\eqdef{t+/t}$$
which of course implies that almost all integers have two divisors satisfying $d_1<d_2<2d_1$. It would be of some interest to get an asymptotic formula for 
$$\T(x):=\sum_{n\leqslant x}\tau^+(n).\eqdef{avt+}$$
It is easy to prove that $\T(x)/(x\log x)\to1$.}
\par \smallskip
This is an example of Erd\H os' way of attacking conjectures from many different angles. Indeed, it is often the case that a stronger statement is more accessible than a weaker one, because it reveals a deeper feature. Here, $\tau^+(n)<\tau(n)$ would suffice to prove the desired conjecture, but Erd\H os asks for much more. As it turns out,  hypothesis \eqref{t+/t} is wrong (and the constant 1 in the last statement should be replaced by 0, most certainly a lapsus digiti), but the idea of considering the measure of the set $\cup_{d|n}(\log d+[-\dm,\dm])$ was precisely that which  eventually led to the solution in \citer{MT84}.\par 
Improving on an estimate of \citer{ET81} that was already sufficient to invalidate \eqref{t+/t}, it was shown in \citer{HT88} (Chapter 4) that the arithmetic function $\tau^+(n)/\tau(n)$ has a limiting distribution $\nu(z)$ satisfying
$${z\over \sqrt{\log (2/z)}}\ll\nu(z)\ll z\log (2/z)\qquad (0<z<1).\eqdef{frt+/t}$$
Thus, $\nu$ is certainly continuous at the origin. Two interesting open problems are (i) to improve upon \eqref{frt+/t} and (ii) to determine, if any, the discontinuity points of the distribution function $\nu$. 
\par
Regarding the second question, I can prove the following.
\Propt{cont1}{The distribution function $\nu$ is continuous at $z=1$.}
\nid We know from theorem 51 of \citer{HT88} (but this already follows from the analysis given in \citer{MT84}) that, for every $\varepsilon>0$, there exists  $T_\varepsilon>\e^{1/\varepsilon}$ such that all integers $n$ except at most those from a sequence of upper density $\leqslant \varepsilon/3$ have two divisors $d$, $d'$, such that $d<d'<2^\varepsilon d<T_\varepsilon$. We may of course assume that $T_\varepsilon$  increases with $1/\varepsilon$. Write $n_\varepsilon:=\prod_{p^j\|n,\,p\leqslant T_\varepsilon}p^j$. For a non-exceptional integer~$n$ and each $m|(n/n_\varepsilon)$, the two divisors $md$ and $md'$  belong to the same interval $]2^k,2^{k+1}[$ $(k\in\N)$ unless $|(\log md)/\log 2-k-1|<\varepsilon$. However, it has been shown in lemma 48.1 of \citer{HT88} that the discrepancy of the sequence $\{(\log m)/\log 2:m|(n/n_\varepsilon)\}$ does not exceed $\varepsilon$ on a subsequence of lower density $1-\varepsilon/3$. Thus, if we discard a sequence of integers $n$ of upper density at most $2\varepsilon/3$, we have
$$\tau^+(n)\leqslant \tau(n)-(1-\varepsilon)\tau(n/n_\varepsilon).$$
Since, for instance, $\tau(n_\varepsilon)\leqslant \log T_\varepsilon$ holds on a sequence of lower density $1-\varepsilon/3$, we get that $$\tau^+(n)\leqslant \tau(n)\Big\{1-{1\over 2\log T_\varepsilon}\Big\}$$
except at most on a sequence of upper density $\varepsilon$. Writing $\eta:=1/\{2\log T_\varepsilon\}$, we have therefore proved that $\nu(1-\eta)\geqslant 1-\varepsilon=\nu(1)-\varepsilon$. Observing that $\varepsilon$ tends to 0 as a function of $\eta$, we obtain the required result.\qed
\par  
According to a copy of the galley-proof that Nicolas forwarded to me at the time, the  statement concerning $\T(x)$ is probably due to some last-minute confusion. It is nevertheless linked to another very interesting problem in probabilistic number theory.
\par 
Let $H(x,y,z)$ denote the number of integers not exceeding $x$ having a divisor in $]y,z]$, so that, with the notation \eqref{avt+},
$$\T(x)=\sum_{2^k\leqslant x}H\big(x,2^k,2^{k+1}\big).$$  There is a large literature on $H(x,y,z)$, starting with \eqref{Bes} and \eqref{E/B}, which can already be seen as evaluations of $$\limsup_{T\to\infty}\lim_{x\to\infty} H(x,T,2T)/x,\hbox{ and }\lim_{T\to\infty}\lim_{x\to\infty}H(x,T,T^{1+\varepsilon_T})/x,$$ respectively. We refer the reader to the recent paper \citer{Fo08} for the history of estimates of $H(x,y,z)$ in the various ranges of the parameters. Here, we only quote the evaluation
$$H(x,y,2y)\asymp {x\over (\log y)^{\delta}(\log_2y)^{3/2}}\qquad (2\leqslant y\leqslant \sqrt{x})\eqdef{FH}$$
with $\delta:=1-(1+\log_22)/\log 2\approx0.08607$. These bounds improve on those of \citer{Te84}, where it is shown by a much simpler analysis that $\e^{-c_1\sqrt{\log _2y}}\leqslant H(x,y,2y)(\log y)^\delta/x\leqslant c_2/\sqrt{\log_2 y}$ for suitable constants $c_1$, $c_2$. Using the symmetry of the divisors of $n$ around $\sqrt{n}$, we easily deduce from \eqref{avt+} and \eqref{FH} the following estimate proved in \citer{Fo08}:
$$\T(x)\asymp {x(\log x)^{1-\delta}\over (\log_2x)^{3/2}}\cdot\eqdef{og-avt+}$$
\par 
Thus, we still fall short of an asymptotic formula for $\T(x)$, although we are now fairly close to one---another challenging problem from an old paper.
\par \smallskip
Let us continue.\par 
{\it Another interesting and unconventional problem states as follows: let $1=d_1<d_2<\cdots<d_{\tau(n)}=n$ be the set of divisors of $n$. Put $$\G(n):=\sum_{1\leqslant i<\tau(n)}{d_i\over d_{i+1}}\cdot$$
I conjecture that $\G(n)\to\infty$ if we disregard a sequence of integers $n$ of density  0. This again would imply the conjecture on $d_1<d_2<2d_1$, but needless to say I cannot prove it. \par 
It would be of interest to determine the normal order of $\tau^+(n)$ and of $\G(n)$ (or at least of $\log \tau^+(n)$ and $\log \G(n)$). Also an asymptotic formula for 
$$\sum_{n\leqslant x}\G(n)$$
would be of interest. It is easy to prove that $(1/x)\sum_{n\leqslant x}\G(n)\to\infty$.}
\par \smallskip
It turns out to be almost trivial that $\G(n)\to\infty$ \pp. Indeed, if $p$ is the smallest prime factor of $n$, then $pd_i|n$ for at least $\dm\tau(n)$ values of $i$ and hence $\G(n)>\tau(n)/2p$. In particular, we have $\G(n)>\tau(n)/\xi(n)$ $\pp$ whenever $\xi(n)\to\infty$. It is, however, not true that this lower bound implies \eqref{conjE}.  Erd\H os  probably had in mind the correct statement that \eqref{conjE} follows from $\G(n)>\dm\tau(n)$ $\pp$, in other words that the distribution function of $\G(n)/\tau(n)$, if it exists, is supported on $[\dm,1]$. 
\par\smallskip 
Erd\H os and I proved in \citer{ET83} that $\G(n)/\tau(n)$ does have a distribution function. We actually established  a fairly general statement: given any bounded real function $\vartheta$ defined on $]0,1[$, the arithmetical function
$$F(n;\vartheta):={1\over \tau(n)}\sum_{1\leqslant i<\tau(n)}\vartheta\Big({d_i\over d_{i+1}}\Big)$$
has a limiting distribution.\note{Note that, in the case $\vartheta:=\1_{[1/2,1]}$, the continuity at 0 of this distribution follows from  \ref{cont1} above and in turn implies \eqref{conjE}. This, however, does not yield a new proof of \eqref{conjE} since we actually used a refinement of \eqref{conjE}  to establish \ref{cont1}.}
\par 
But it is not true that the distribution function of $\G(n)/\tau(n)$ is supported on $[\dm,1]$. Indeed, we can show that 
$$\dd\{n\geqslant 1:\G(n)/\tau(n)\leqslant \varepsilon\}>0\qquad (0<\varepsilon\leqslant 1).$$
This follows from the fact that most integers $n$ free of small prime factors are such that $d_i<\dm\varepsilon d_{i+1}$ for most indices $i$. We omit the details, which  can easily be reconstructed from lemma 4 of \citer{ET81} and lemma 3 of \citer{ET83}. \par 
As far as average orders are concerned, it is proved in \citer{ET83} that 
$$\sum_{n\leqslant x}F(n;\vartheta)=\vartheta(1)x\log x+O\Big({x(\log x)^{1-\delta}\log_3x\over \sqrt{\log_2x}}\Big),$$
provided $\vartheta$ is twice continuously differentiable on $[0,1]$. Here $\delta$ is as in \eqref{FH} and the exponent of $\log x$ is optimal. Moreover, by theorem 3 of \citer{ET83} and \eqref{FH}, we obtain the improvement
$${c_1x(\log x)^{1-\delta}\over (\log_2x)^{3/2}}\leqslant x\log x-\sum_{n\leqslant x}\G(n)\leqslant {c_2x(\log x)^{1-\delta}\over (\log_2x)^{3/2}},$$
valid for suitable positive constants $c_1$, $c_2$.
\smallskip
After a discussion on the normal size of the $k$-th prime factor $p_k(n)$ of an integer $n$ and a simple proof, via the Tur‡n--Kubilius inequality, of the asymptotic formula
$$\log_2p_k(n)\sim k\qquad (k\to\infty)\qquad \pp,\note{We do not reproduce this and  refer the reader to \citer{HT88} (chapter 1) and to \citer{Te08} (theorem III.3.10).}\eqdef{log2pk}$$ Erd\H os describes a problem on fractional parts of Bernoulli numbers, which does not fit with the focus of this survey. Then, he states two problems related to densities of integer sequences.
\smallskip\goodbreak
{\it Denote by $\lambda_k(p)$ the density of the integers $n$ whose $k$-th prime factor is $p$. $\lambda_k(p)$ can easily be calculated by the exclusion-inclusion principle (essentially the sieve of Eratosthenes). By \eqref{log2pk}, for almost all integers, $p_k(n)$ is about $\exp\exp k$. On the other hand, it is easy to see that the largest value of $\lambda_k(p)$ is assumed for much smaller values of $p$, in fact for
$$\e^{k(1-\varepsilon)}<p<\e^{k(1+\varepsilon)}.$$
By more careful computation it would easily be possible to obtain better estimates. The simple explanation for this apparent paradox is that there are very many more values of $p$ at $\e^{\e^k}$ than at $\e^k$. It is not impossible that $\lambda_k(p)$ is unimodal, \ie it first increases with $p$, then assumes its maximum and then decreases. I in fact doubt that $\lambda_k(p)$ behaves so regularly but have not disproved it.\par 
The same problems arise if $\Lambda_k(d)$ denotes the density of the integers $m$ whose $k$-th divisor is $d$. Here I obtain that if $d_1(n)<d_2(n)<\cdots$ are the consecutive divisors of $n$ then for all but $\varepsilon x$ integers $n\leqslant x$ for $k>k_0(\varepsilon,n)$
$$\exp\big\{k^{(1/\log 2)-\varepsilon}\big\}<d_k(n)<\exp\big\{k^{(1/\log 2)+\varepsilon}\big\}.$$
On the other hand, for fixed $k$, $\Lambda_k(d)$ is maximal for
$$\e^{(1-\varepsilon)\log k\log_2k}<d<\e^{(1+\varepsilon)\log k\log_2k}.\eqdef{EmodeL}$$
It can be shown that $\Lambda_k(d)$ is not unimodal.}
\smallskip
The existence of the densities $\lambda_k(p)$ and $\Lambda_k(d)$ immediately follows from the fact that the sequences under consideration are finite unions of congruence classes. The idea of considering the local laws of the distributions of $p_k(n)$ and $d_k(n)$ stems naturally from the law of iterated logarithm underlying \eqref{log2pk} (and based upon the fact that the variables $U_j(n)$ defined in~\eqref{Ujn} are almost Gaussian): indeed, Erd\H os announced in 1969 \citer{Er69} that
$$\sum_{\log_2p\leqslant k+z\sqrt{k}}\lambda_k(p)=\Phi(z)+o(1)\qquad (k\to\infty),\qquad \Phi(z):={1\over \sqrt{2\pi}}\int_{-\infty}^z\e^{-t^2/2}\d t.\eqdef{Esplk}$$
Thus, the study of the $\lambda_k(p)$ is another way of looking at the asymptotic independence of the small prime factors, while, as it turns out, the study of the $\Lambda_k(d)$ is a (positive) test of the dependence of the divisors.
\par
By the sieve of Eratosthenes, we have
$$\lambda_k(p)={1\over p}\prod_{q<p}\Big(1-{1\over q}\Big)s_{k-1}(p)\qquad (k\geqslant 1),\eqdef{lamkn}$$
where $q$ denotes a prime number and we have put 
$$s_j(p):=\sum_{\di{P^+(m)<p}{\omega(m)=j}}{1\over m}\quad(j\geqslant 0).$$
 Thus, we have identically
$$F(z,p):=\prod_{q<p}\Big(1+{z\over q-1}\Big)=\sum_{j\geqslant 0}s_j(p)z^j.$$
As noted by Balazard,\note{Private communication, February 28, 1989.} this settles, in the affirmative, the question of the unimodality of the sequence $\{s_j(p)\}_{j\geqslant 1}$ and hence of $\{\lambda_k(p)\}_{k}$ for all $p$. Indeed, it is well known (see, \eg, \citer{PS98}, Part V, problem 47) that, if a polynomial has only real roots, then the number of sign changes in the sequence of its coefficients is equal to the number of positive roots. Since, for all positive numbers $a_1,\ldots, a_n$, the polynomial 
$$(1-x)\prod_{1\leqslant j\leqslant n}(x+a_j)=\sum_{0\leqslant r\leqslant n+1}(\sigma_{n-r}-\sigma_{n+1-r})x^r$$
where $\sigma_h:=\sum_{1\leqslant j_1<j_2<\cdots<j_h\leqslant n}a_{j_1}\cdots a_{j_h}$ $(0\leqslant j\leqslant n+1)$, has exactly one positive root, it follows that the sequence $\{\sigma_h\}_{h=0}^{n}$ of elementary symmetric functions of the $a_j$ is unimodal. Applying this with $\{a_j\}_{j=1}^{n}=\{1/(q-1):q<p\}$ yields the stated property.
\par 
Of course the above argument tells us nothing about the mode. An analysis of $\lambda_k(p)$ by the saddle-point method has been achieved by Erd\H os and myself in \citer{ET89a}. I only quote a few results from this work. Write
$$L:=\log \Big({\log p\over \log (k+1)}\Big), \quad M:=\log\Big({\log p\over 1+\log ^+(k/L)}\Big),\quad R:=L\big\{1+\log ^+(k/L)\big\}.$$
Then, given any $\varepsilon>0$, we have
$$\eqalign{\lambda_k(p)&={1\over p}\prod_{q<p}\Big(1-{1\over q}\Big){M^{k-1}\over (k-1)!}\e^{O((k-1)/R)}\cr
{\lambda_{k+1}(p)\over \lambda_k(p)}&={M\over k}\Big\{1+O\Big({M\over R}\Big)\Big\}\cr}\qquad (1\leqslant k\leqslant p^{1-\varepsilon}).$$
Moreover, we have, for all primes $p$,
$$\max_{k\geqslant 1}\lambda_k(p)={1+O(1/\log_2p)\over p\sqrt{2\pi\log_2p}}$$
and any value of $k$ realizing this maximum satisfies $k=\log_2p+O(1)$. \par 
For fixed $k$, the result we found was slightly different from that foreseen by Erd\H os, probably through a hasty computation. We actually have
$$\max_p\lambda_k(p)=\exp\Big\{-k\Big(\log k-\log_2k-1+{2\log_2k+1\over \log k}+{2(\log_2k)^2-\log_2k+O(1)\over (\log k)^2}\Big)\Big\}$$
and any value of $p $ realizing  this maximum satisfies
$$\log p={k\over \log k}\Big\{1+{2\log_2k\over \log k}+{2(\log_2k)^2-3\log_2k+O(1)\over (\log k)^2}\Big\}.$$
It remains that the phenomenon described by Erd\H os does hold: modal values of the sequence $\{\lambda_k(p)\}_p$ occur at relatively small values. In other words, in the series
$$\sum_{p}\lambda_k(p)=1$$
the decrease of the general term as a function of $p$ is so slow that the contribution of the very numerous terms around $\exp\exp k$ dominate, while the `large' values around $\e^{k/\log k}$ are too few, and indeed not sufficiently large, to contribute significantly to the sum.
\par 
To my knowledge, the problem of the (probably non) unimodality of the sequence $\{\lambda_k(p)\}_p$ is still open.
\par \goodbreak
In \citer{DKT02}, De Koninck and I improve on \eqref{Esplk}. Uniformly for $k\geqslant 1$, $z\in\r$, we have
$$\sum_{\log_2p\leqslant k+z\sqrt{k}}\lambda_k(p)=\Phi(z)+{\Phi_0(z)\over \sqrt{2\pi k}}+O\Big({1\over k}\Big)$$
with $$\Phi_0(z):=\e^{-z^2/2}\{\ft13+A-\ft13z^2\},\quad A:=\gamma-\sum_{p}\Big\{\log \Big({1\over 1-1/p}\Big)-{1\over p}\Big\}\approx 0.26150.$$
Here $\gamma$ denotes Euler's constant.\par 
This yields estimates for the median value of the distribution of the $k$-th prime factor, defined as the largest prime $p^*=p^*_k$ such that $\sum_{p\leqslant p^*_k}\lambda_k(p)< \dm$. We find that
$$\log_2p_k^*=k-b+O(1/k)\qquad (k\geqslant 1)\eqdef{medpk}$$
with $b=\ft13+A\approx0.59483$. Numerical computations provide $p^*_2=37$, $p_3^*=42719$.
\par \smallskip
A clear descendant of this problem  is the following formula, also proved in \citer{DKT02}, which turns out to be an application of the estimate for partial sums of the exponential series---an ancient problem of Ramanujan---needed to prove \eqref{medpk}. We have
$$\sum_{\di{n\leqslant x}{\Omega(n)\leqslant \log_2x}}1=\dm x-x{C+\chv{\log_2x}\over \sqrt{2\pi\log_2x}}+O\Big({1\over \log_2x}\Big)\qquad (x\geqslant 3),$$
where 
$C:=A-\ft23-\sum_{p}1/\{p(p-1)\}\approx0.36798$ and $\chv t$ denotes the fractional part of the real number $t$.
\par 
As is to be expected, the results on $\Lambda_k(d)$ are much less precise. Erd\H os' pp-estimate for $\{d_k(n)\}_{1\leqslant k\leqslant \tau(n)}$ immediately follows from the law of iterated logarithm for the prime factors. We obtain in particular, for all $\varepsilon>0$,
$$\sum_{|\log_2d-(\log k)/\log 2|>R_k}\Lambda_k(d)=o(1)\qquad (k\to\infty),$$
with $R_k:=\sqrt{\{(2+\varepsilon)/\log 2\}\log k\log_3k}$. Thus, we can consider that the problem of normal order of $d_k(n)$ is essentially solved. In \eqref{EmodeL}, Erd\H os raises the problem of modal values of $\Lambda_k(d)$ \ie of determining as precisely as possible those $d$ such that $$\Lambda_k(d)=\Lambda_k^*:=\max_m\Lambda_k(m).$$ He announces a result which we shall see to be slightly incorrect but nevertheless unveils a rather deep phenomenon.
\par 
Let $\tau(n,z)$ denote the number of divisors of $n$ not exceeding $z$. The following formula, proved in \citer{ET89a}, is the analogue of \eqref{lamkn}:
$$\Lambda_k(d)={1\over d}\prod_{p\leqslant d}\Big(1-{1\over p}\Big)\sum_{\di{P^+(m)\leqslant d}{\tau(md,d)=k}}{1\over m}\cdot$$
Here, the $m$-sum obviously depends on the arithmetic structure of $m$ and seemingly harmless questions may reveal to be quite delicate, such as the proof given in \citer{ET89a} of the equivalence
$$\Lambda_k(d)>0\iff\tau(d)\leqslant k\leqslant d.\eqdef{Lkd=0}$$
Let us put
$$K_j:=k^{(\log_{j+2}k)/\log 2}\qquad (j\geqslant 0).$$
It is well known that $\min_{\tau(d)\geqslant k}d=K_0^{1+o(1)}.$
Now let $N_y:=\prod_{p\leqslant y}p$, where $y$ is the smallest integer such that $\tau(N_y)=2^{\pi(y)}\geqslant k$. By selecting $d=d_k(N_y)$ and reducing the $m$-sum above to the single value $m=N_y/d$, we obtain the left-hand side of the double inequality
$${k^{O(1)}\over K_0K_1}\leqslant \Lambda_k^*\leqslant {k^{O(1)}K_1\over K_0}$$
proved in \citer{ET89a}, while the upper bound already needs a rather involved analysis of the sum. This led Erd\H os and I in \citer{ET89a} to express the belief that the correct version of \eqref{EmodeL} should be
$d=K_0^{1+o(1)}$.
\par 
Indeed, there are essentially two sound models for the structure of those $d$ realizing the mode. Either $\tau(d)\approx k$ and hence $d\approx K_0$ and therefore the $m$-sum has size $\asymp1$, or $m$ and $d$ contribute evenly to the divisors counted by $\tau(md,d)$ and $\tau(d)\approx \tau(m,d)\approx \sqrt{k}$, so that $d$ and the values of $m$ appearing in the sum are all at least of size $\sqrt{K_0}$. This latter possibility is of course much more complex than the former, since it implies the existence of many integers $m$ having divisors combining with those of $d$ in such a way that $\tau(md,d)=d$. The above belief corresponded to the conviction that the simplest situation did prevail. However, in \citer{BT02}, La Bretche and I show that this is not the case: for large $k$, we have
$${k^{O(1)}\over K_0\sqrt{K_1}K_2}\leqslant \Lambda_k^*\leqslant {\sqrt{K_2}k^{O(1)}\over K_0\sqrt{K_1}},\qquad \Lambda_k(d)=\Lambda_k^*\Rightarrow d=K_0^{1/2+o(1)}.$$
(See \citer{BT02} for a more precise statement and some further information.)\par 
Here again, Erd\H os' question led to a deeper understanding of the structure of the set of divisors of certain classes of integers and revealed an unexpected phenomenon.
\par 
The conjecture \eqref{EmodeL}, although inaccurate, clearly satisfies all criteria quoted at the beginning of this paper. As far as criterion (iv) is concerned, we quote from \citer{BT02} the following estimate, where $\Psi_1(x,y)$ denotes the number of $y$-friable squarefree integers not exceeding~$x$, \ie $\Psi_1(x,y):=\sum_{n\leqslant x,\,P^+(n)\leqslant y}\mu(n)^2$. Given any $\kappa\geqslant 1$, we have
$$\Psi_1(x+x/z,y)-\Psi_1(x,y)\ll\Psi_1(x,y)/z\qquad \big(x\geqslant 2,\,y\geqslant 2,\,1\leqslant z\leqslant \min(x,y^\kappa)\big).$$
\par 
The statement concerning the non-unimodality of $\{\Lambda_k(d)\}_{d}$ follows easily from \eqref{Lkd=0}, since, for any $\varepsilon>0$, we can construct four integers such that $$K_0^{1+\varepsilon}<p_1<d_1<p_2<d_2<2K_0^{1+\varepsilon},$$ where the $p_j$ are primes and the $d_j$ satisfy $\tau(d_j)>k$ and hence $\Lambda_k(d_j)=0$ $(j=1,2)$.
\medskip
In the next paragraphs of \citer{Er79}, Erd\H os quotes a number of results  related to the normal distribution of prime factors, some of which are stated in \citer{Er69}. For instance, he explains that, with the notation \eqref{Ujn}, the statement that $U_j(n)$ and $U_h(n)$ are asymptotically independent provided $j/h\to\infty$ follows from the methods of \citer{EK40}, his epoch-making paper with Kac on the Gaussian distribution of prime factors.  He also comments on the fact that \eqref{log2pk} shouldn't be taken too literally by stating the following theorem, which I reproduce with a few changes in the notation.\par \smallskip
{\it Let $\{\alpha_k\}_{k=0}^{\infty}$ tend monotonically to 0 as $k\to\infty$. Denote by $h_\alpha(n)$ the number of $k$ such that $|\log_2p_k(n)-k|\leqslant \alpha_k$. Then, if $\sum_{k}\alpha_k/\sqrt{k}<\infty$, for every integer $m$ the set $\{n\geqslant 1:h_\alpha(n)=m\}$ has a natural density $\beta_m$ and $\sum_{m}\beta_m=1$, in other words $h_\alpha(n)$ has a limiting distribution, while, if $\sum_{k}\alpha_k/\sqrt{k}=\infty$, $h_\alpha(n)\to\infty$ \pp.}\smallskip\par\goodbreak 
As far as I know, none of these results has ever been proved in full detail and no effective versions of the statements have been investigated. It would be quite interesting to pursue these tasks with the powerful analytical tools that have been devised since Erd\H os' paper was written.
\medskip
The next section of \citer{Er79} introduces a fundamental concept.
\smallskip
{\it Let $p_1<p_2<\ldots$ be an infinite sequence of primes. It is quite easy to prove that
$$\sum_{}{1\over p_i}=\infty$$
is the necessary and sufficient condition that almost all integers $n$ should have a prime factor $p_i$. It seems very difficult to obtain a necessary and sufficient condition that if $a_1<a_2<\ldots$ is a sequence of integers then almost all integers $n$ should be a multiple of one of the $a$'s. \par 
I just want to illustrate the difficulty by a simple example. Let $n_{i+1}>(1+c)n_i$. Consider the integers $m$ which have a divisor $d$ satisfying $n_k<d\leqslant n_k(1+\eta_k)$. If $\sum_{k\geqslant 1}\eta_k<\infty$, then it is easy to see that the density of these integers exists and is less than 1. If $\sum_{k\geqslant 1}\eta_k=\infty$, it seems difficult to get a general result, \eg\ if $\eta_k=1/k$ the density in question exists and is less than 1. It seems certain that there is an $\alpha$, $0<\alpha<1$, so that if $\beta<\alpha$ and $\eta_k=1/k^\beta$ the density of the $m$ having a divisor $d$, $n_k<d\leqslant n_k(1+\eta_k)$ is 1 and if $\beta>\alpha$ it is less than 1.}
\smallskip
The problem raised here may be reformulated as follows: characterise those integer sequences $\A$ such that $\d\M(\A)=1$. Following Hall \citer{Ha90}, we call such a sequence $\A$ a {\it Behrend sequence}. This concept has been a constant concern for Erd\H os during more than fifty years: while, as he remarks in the above excerpt, the corresponding problem is easy when one considers a sequence of primes, or, more generally, a sequence of pairwise coprime integers, delicate and interesting questions arise immediately in the general case, corresponding to the study of strongly dependent random variables.\par 
By the Davenport--Erd\H os theorem \citer{DE37} quoted earlier, a necessary and sufficient condition that $\A$ should be a Behrend sequence is that $\delta\M(\A)=1$ where $\delta$ stands for the logarithmic density. Thus, we have obviously that $\delta \A=1$ is a sufficient condition for $\A$ to be a Behrend sequence. For a long time, I thought that this should have a simple, direct proof, but I could not find one that wasn't essentially equivalent to the Davenport--Erd\H os general and deep result that $\underline\dd\M(\A)=\delta\M(\A)$ for any sequence $\A$.
I eventually came up with the following.
\Propt{DEbase}{Let $\A$ be an integer sequence such that $\delta\A=1$. Then $\d\M(\A)=1$.}
\nid
Recall that we defined $P^+(n)$ as the largest prime factor of an integer $n$ with the convention that $P^+(1)=1$. Symmetrically we let $P^-(n)$ denote the smallest prime factor of $n$ and set $P^-(1)=\infty$. 
For $y\geqslant 1$, let us write $\A_y:=\{n\in\A:P^+(n)\le y\}$, $n_y:=\prod_{\di{p^\nu\|n}{p\le y}}p^\nu$.
As
$n_y\in\M(\A)$   implies $n\in\M(\A)$, we plainly have for any fixed $y\ge 1$ and $x\to\infty$,
$$\eqalign{{1\over  x}\sum_{\di{n\le x} {n\in\M(\A)}}1\ge {1\over
x}\sum_{\di{r\in\M(\A_y)}{P^+(r)\le y}}\sum_{\di{s\le
x/r}{P^-(s)>y}}1&={1\over
x}\sum_{\di{r\in\M(\A_y)}{P^+(r)\le y}}\Big\{{x\over r}\prod_{p\le y}\Big(1-{1\over
p}\Big)+O(1)\Big\}\cr
&\to m(y):=\prod_{p\le y}\Big(1-{1\over
p}\Big)\sum_{\di{r\in\M(\A_y)}{P^+(r)\le y}}{1\over r}\cdot\cr}
$$
Thus, we only have to show that $m(y)\to1$ as $y\to\infty$. We have trivially
$$m(y)\ge \prod_{p\le y}\Big(1-{1\over
p}\Big)\sum_{r\in\A_y}{1\over r}\cdot$$
Writing $a$ for an element of $\A$, we deduce from our hypothesis $\delta\A=1$ that there is a non-increasing function
$\epsilon(x)$ tending to 0 as
$x\to\infty$ such that, for  $1\le y\le x$,
$$\{1-\epsilon(x)\}\log x\le \sum_{a\le x}{1\over a}\le \sum_{r\in \A_y}{1\over
r}+\sum_{\di{n\le x}{P^+(n)>y}}{1\over n}\cdot\eqdef{1}
$$
Setting $u:=(\log x)/\log y$, we may rewrite the last sum in \eqref{1} as
$$\eqalign{\sum_{n\le x}&{1\over n}-\sum_{P^+(n)\le y}{1\over
n}+\sum_{\di{n>x}{P^+(n)\le y}}{1\over n}\cr&\le \log x+O(1)-\prod_{p\le
y}\Big(1-{1\over p}\Big)^{-1}+x^{-1/\log y}\prod_{p\le y}\Big(1-{1\over
p^{1-1/\log y}}\Big)^{-1}\cr&\le \big\{1+O(\e^{-u}\big\}\log x+O(1)-\prod_{p\le
y}\Big(1-{1\over p}\Big)^{-1}.\cr}
$$
Inserting back into \eqref{1}, we get
$$\sum_{r\in \A_y}{1\over
r}\ge \prod_{p\le
y}\Big(1-{1\over p}\Big)^{-1}+O\big(1+\{\e^{-u}+\epsilon(x)\}u\log y\big) . $$
It remains to select $u=1/\sqrt{\epsilon(y)}$ and let $y\to\infty$ to obtain $\lim_ym(y)=1$.
\qed
The paper \citer{DE51} (see also \citer{Te08}, Exercises 247-249) contains another fundamental formula,~viz.
$$\un\dd\M(A)=\lim_{T\to\infty}\dd\M(\A\cap[1,T]).\eqdef{dseq}$$ 
We call the right-hand side the {\it sequential density} of the set of multiples $\M(\A)$. From Behrend's fundamental inequality, valid for finite sequences, we hence deduce from \eqref{dseq} that
$$1-\un\dd\M(\A\cup \B)\geqslant \{1-\un\dd\M(\A)\}\{1-\un\dd\M(\B)\}\eqdef{Beh-in}$$
holds for all integer sequences $\A$, $\B$.\note{This has been nicely improved by Ahlswede and Khachatrian \citer{AK95}.} It follows in particular that
$$\sum_{a\in\A}{1\over a}=\infty\eqdef{CNB}$$
is a necessary condition for $\A$ to be a Behrend sequence, and that any tail $\A\sset[1,T]$ of a Behrend sequence is still a Behrend sequence.
\par 
The structure of Behrend sequences long intrigued Erd\H os. The problem is indeed quite intricate and even seemingly innocent questions, such as that of a criterion for $\A$ to be a Behrend sequence in the special case when the members of $\A$ only have a bounded number of, or even at most two, prime  factors, do not have a simple answer: such a criterion is given in Ruzsa--Tenenbaum \citer{RT96} in the case of two prime factors; in Erd\H os--Hall--Tenenbaum \citer{EHT94}, it is shown that $\dd\M(\A)$ always exists when the number of prime factors is bounded but that this condition is optimal.
\par 
Another interesting feature of Behrend sequences, proved in \citer{HT92}, is that, if $\A$ is a Behrend sequence, then $\sum_{d|n,\,d\in\A}1\to\infty$ \pp.\par \smallskip
Since it seems hopeless to find an effective criterion for the general situation, we are led to consider sequences with a special structure. The sequence $\E$ in \eqref{dd'} is one example. Another instance is that of block sequences, appearing implicitly in Erd\H os' formulation above. As in \citer{HT92}, we formally define a block sequence by the property that it can be written in the form
$$\A=\bigcup_{j\geqslant 1}\A_j,\qquad \A_j:=]T_j,H_jT_j]\cap\N^*\quad(j\geqslant 1),$$
where the (disjoint) blocks $\A_j$ satisfy some growth condition that guarantees some local regularity, namely that, for some fixed parameter $\eta>0$,
$$1+1/T_j^{1-\eta}\leqslant H_j\leqslant \min(T_j,T_{j+1}/T_j)\quad(j\geqslant 1).\eqdef{def-bs}$$
When the $T_j$ grow sufficiently fast, we might then expect a Borel--Cantelli type criterion enabling us to decide whether $\A$ is a Behrend sequence according to whether a certain series involving the quantities $\dd\M(\A_j)$ diverges or not.\par 
These questions have had a fairly wide posterity and many descendants. We refer the reader, in particular, to the papers \citer{Ha90}, \citer{HT92}, \citer{RT96}, \citer{Te96a} and to the book \citer{Ha96}, for a number of results and conjectures on Behrend sequences and uniform distribution on divisors. Here, we only quote two significant results which confirm, at least in the case of block sequences, that a criterion of Borel--Cantelli type  is relevant.
\par 
In order to avoid technical hypotheses, we restrict to special cases which still reflect the general picture. We start with a result concerning the situation when the blocks are somewhat short.\note{See \citer{Te96} for an explanation of the fact that any criterion for block Behrend sequences can be split into one in which the block are assumed to be short, in some precise way, and one in which the blocks are assumed to be long.} The necessity part is due to Hall--Tenenbaum \citer{HT92}, and the sufficiency to Tenenbaum \citer{Te96}.
\propt{CNSB}{ (\citer{HT92}, \citer{Te96})}{Let $\A=\cup_j\A_j$ be a block sequence such that, for suitable real constants $\alpha$, $\gamma$, $\sigma$, $\tau$, with $\sigma>-1$, we have 
$$\log (T_{j+1}/T_j)\asymp j^\sigma(\log j)^{\tau},\qquad \log H_j\asymp (\log j)^\gamma/j^\alpha\qquad (j\to\infty).$$
Put $\sigma_0:=(\log 2)/(1-\log 2)$ and define  $$\alpha_0(\sigma):=\cases{(1-\log 2)(\sigma_0-\sigma)&if $-1<\sigma\leqslant \sigma_0$,\cr
\sigma_0-\sigma& if $\sigma>\sigma_0$.\cr}$$
Then $\A$ is a Behrend sequence if $\alpha<\alpha_0(\sigma)$ and $\A$ is not a Behrend sequence if $\alpha>\alpha_0(\sigma)$.}
Note that \eqref{def-bs} implies $\sigma+\alpha>0$ or $\sigma+\alpha=0$ and $\gamma\leqslant \tau$.
\par If we set $\sigma=\tau=\gamma=0$, we obtain that, provided $1+c_1\leqslant T_{j+1}/T_j\leqslant 1+c_2$ for suitable constants $c_1>0$, $c_2>0$, and $H_j:=1+1/j^\alpha$ $(j\geqslant 1)$, then the block sequence $\A$ is a Behrend sequence if $\alpha<\log 2$ and is not if $\alpha>\log 2$. This settles Erd\H os' conjecture quoted above. His original claim was that the critical exponent $\alpha_0$ should exist under the sole condition $T_{j+1}/T_j>1+c_1$, but this cannot hold as it stands since it follows from theorem 1 of \citer{HT92} that $\A$ is not a Behrend sequence for any $\alpha$ if we set, for instance, $T_j:=\exp\exp j$ $(j\geqslant 1)$. However, he explained later on, in private conversation, that he really had in mind a two-sided condition. 
\par 
From Behrend's inequality \eqref{Beh-in}, the condition
$$\sum_{j}\dd\M(\A_j)=\infty$$
is necessary for a block sequence to be a Behrend sequence. However, this is in general much weaker than the sufficient condition obtained in \citer{HT92}. For instance, if we assume, in the setting of \ref{CNSB}, that $-\sigma<\alpha\leqslant 0$ or that $\alpha=-\sigma\leqslant 0$ and $\gamma<\tau$, then we have from Ford's estimates in \citer{Fo08} that
$$\dd\M(\A_j)\asymp {(\log 2j)^{(\gamma-\tau)\delta-3/2}\over j^{(\sigma+\alpha+1)\delta}}\qquad (j\geqslant 1),$$
where $\delta$ is  as in \eqref{FH}, while \ref{CNSB} tells us that $\A$ is a Behrend sequence if
$$\sum_{j}{1\over j^{(\sigma+\alpha+1)\beta}}=\infty$$
for some $\beta>1-\log 2$ and, moreover, that $\A$ is not a Behrend sequence if the above series converges for some $\beta<1-\log 2$.
Hence, we have a pseudo Borel--Cantelli criterion of the shape 
$$\sum_{j}\{\dd\M(\A_j)\}^{c+o(1)}=\infty,$$
with $c:=(1-\log 2)/\delta\approx 3.566509$.
It would be very interesting to have a probabilistic interpretation for conditions of this type.
\par 
For the special sequence
$$\A_\lambda:=\bigcup_{j\geqslant 1}]\exp j^\lambda,2\exp j^\lambda]\cap\N^*,$$
a refined approach of the same technique yields in \citer{HT92} a complete proof of Erd\H os' so called $\B_\lambda$-conjecture\note{The name of the conjecture comes from the former notation $\B(\lambda)=\M(\A_\lambda)$.} dating at least from the seventies and referred to in \citer{HT88} pp.\ 49 and 63: $\A_\lambda$ is a Behrend sequence if, and only if, $\lambda\leqslant 1/(1-\log 2)$. This is heuristically justified by the assumption that, for almost all $n$, the numbers $(\log d)^{1/\lambda}$ are uniformly distributed modulo 1 when $d$ runs through the divisors of $n$.\note{This is actually proved in \citer{Te82}. See also \citer{HT86} and \citer{Te96a}.} However, the limiting case $\lambda=1/(1-\log 2)$ is not covered by this argument and indeed needs a more delicate proof.
 \par 
In the same spirit, and as a clear descendant of this class of problems, I quote the theorem of Kerner and myself \citer{KT04}, according to which
$$
\min_{d\mid n}\left\Vert d\vartheta \right\Vert =1/\tau
(n)^{1+o(1)}\quad \pp, \eqdef{KT}
$$
provided the sequence of convergents  $\{p_j/q_j\}_{j=0}^{\infty}$ of the real number $\vartheta$  satisfies $$\log q_{j+1}<(\log q_j)^{1+o(1)}.\eqdef{cond-erthd}$$ Here we used the standard notation $\no t=\min_{n\in\z}|t-n|$. Note that, as explained in~\citer{KT04}, it is easy to construct real numbers $\vartheta$ contravening \eqref{KT}. A challenging open question is to determine precisely the set of real numbers $\vartheta$ such that \eqref{KT} holds. We know from \citer{KT04} that \eqref{cond-erthd} cannot be replaced by $\log q_{j+1}<q_j^{(1-\varepsilon)/\log 2}$ with some $\varepsilon>0$.
\par 
When the blocks are long, in a suitable sense, we obtain a similar pseudo-criterion, but with $c=1$--- hence closer to a classical probabilistic approach.
\propt{CNSB2}{ (\citer{HT92})}{Let $\A$ be a block sequence. Assume that, for some $\varepsilon>0$, we have
$$\log H_{j+1}>2(\log T_{j+1})^{\varepsilon}(\log T_j)^{1-\varepsilon}\qquad (j\geqslant 1).$$
Then $\dsp\sum_{j}\Big({\log H_j\over \log T_j}\Big)^{\delta_1}=\infty$ for some $\delta_1>\delta$
implies that $\A$ is a Behrend sequence, while
$\dsp\sum_{j}\Big({\log H_j\over \log T_j}\Big)^{\delta_2}<\infty$ for some $\delta_2<\delta$
implies that $\A$ is not a Behrend sequence.}
We refer the reader to chapter 1 of \citer{Ha96} for further results and comments on Behrend sequences. Once more, we see how fertile Erd\H os' problems and conjectures revealed themselves along the years.
\medskip\goodbreak
Erd\H os follows with refined questions concerning the set of multiples of an interval. I slightly alter the notation in order to match subsequent works.
\medskip
{\it Denote by $\varepsilon(y,z)$ the density of integers having a divisor $d$ satisfying $y<d\leqslant z$ and by $\varepsilon_1(y,z)$ the density of integers having precisely one divisor $d$, $y<d\leqslant z$. Besicovitch proved $\liminf\varepsilon(y,2y)=0$ and I proved that if $(\log z)/\log y\to1$, then $\lim\varepsilon(y,z)=0$ \citer{HR74} (chapter V). It is easy to see that this result is best possible, \ie\ $\lim \varepsilon(y,z)=0$ implies $(\log z)/\log y\to1$.\par 
Further I can prove that $\varepsilon_1(y,z)<c/(\log y)^\alpha$ for a certain $0<\alpha<1$. Perhaps $\varepsilon_1(y,z)$ is unimodal for $z>y+1$, but I know nothing about this. I do not know where $\varepsilon_1(y,z)$ assumes its maximum.\par 
I am sure that $\varepsilon_1(y,z)/\varepsilon(y,z)\to0$ for $z=2y$. If $z-y$ is small then clearly $\varepsilon_1(y,z)/\varepsilon(y,z)\to1$ and I do not know where the transition occurs.
\par 
Some time ago the following question occurred to me: let $k$ be given and $n>n_0(k)$. Is there an absolute constant $\alpha$ so that for every $n<m\leqslant n^k$ there is a $t$, $0<t\leqslant (\log n)^\alpha$, so that $m+t$ has a divisor in $]n,2n]$? More generally: if $n+1=a_1<a_2<\ldots$ is the sequence of integers which have a divisor $d$, $n<d\leqslant 2n$, determine or estimate $\max_{a_i<n^k}(a_{i+1}-a_i)$.\par }
\smallskip
Nearly all these questions are now essentially settled. In \citer{Te87}, I proved that if $z-y\to\infty$ and $z\leqslant y\Big\{1+(\log y)^{1-\log 4}\e^{-\xi\sqrt{\log_2y}}\Big\}$ with $\xi\to\infty$, then $\varrho_1(y,z):=\varepsilon_1(y,z)/\varepsilon(y,z)\to1$, while $\varrho_1(y,z)\geqslant \e^{-c\sqrt{\log y\log_2y}}$ when $z_0(y):=y\big\{1+(\log y)^{1-\log 4}\big\}<z\leqslant 2y$. On seeing this, Erd\H os changed his mind concerning the asymptotic behaviour of $\varrho_1(y,z)$ and conjectured that this quantity should tend to a positive limit for $z=2y$. Ford \citer{Fo08} then proved that $\varrho_1(y,z)\asymp1$ when $y+1\leqslant z\ll y$. Thus,  the transition imagined by Erd\H os should ideally be seen as a frontier between the cases when $\varrho_1(y,z)$ tends to 1 or to a constant less than~1. We still do not know whether $\varrho_1(y,z)$ tends to a limit when $z_0(y)<z\ll y$ but it follows from Ford's estimates in \citer{Fo08} that $\varrho_1(y,z)\to0$ if $z/y\to\infty$. I conjecture that $\varrho_1(y,z)\not\to1$ when $y$, $z$ tend to infinity in such a way that $z>y\big\{1+(\log y)^{1-\log 4+\varepsilon}\big\}$.
\par 
To my knowledge, the question of unimodality of $\varepsilon_1(y,z)$ as a function of $z$ is still open.
\par 
The last problem seems difficult and represents a deep open question. Let $M_n(x)$ denote the counting function of $\M(]n,2n])$ and set 
$$M_n(x)=\varepsilon_nx+R_n(x)\qquad (x\geqslant 1).$$
Then $a_{i+1}-a_i=\{1-R_n(a_{i+1})+R_n(a_i)\}/\varepsilon_n$. Since $1/\varepsilon_n\asymp(\log n)^{\delta}(\log_2n)^{3/2}$, the first question amounts to asking whether
$\max_{a_i\leqslant n^k}|R_n(a_{i+1})-R_n(a_i)|\ll_k(\log n)^\beta$ for some $\beta$ independent of $k$. \par 
Note that Hall \citer{Ha96} studied the quadratic mean of $R_n(x)$. His lower bound implies that $\sup_x|R_n(x)|\gg n^{c}$ with $c:=\ft12-\log (\pi^2/6)/\log 4\approx0.14098$. However, he recently observed \citer{Ha12} that the results obtained in \citer{Ha96} imply much more, namely $$\sup_x|R_n(x)|>2^{\{1+o(1)\}n/(2\log n)}.$$ This follows on noticing that $]n,2n]=\A\cup \B$ where $\A$ comprises all primes in the interval and $\B$ includes all remaining, composite integers. Then $(a,b)=1$ whenever $a\in\A$, $b\in\B$. It only remains to apply equations  (3.26), (3.10) and (3.20) from \citer{Ha96}.\note{The author takes pleasure in thanking Richard R. Hall for letting him include this proof here.} Although this does not contradict Erd\H os' conjecture, it shows that it must be delicate.
\medskip\goodbreak
I conclude this survey of posterity and descendants of Erd\H os' paper \citer{Er79} by quoting a problem that was for him a constant concern even though he thought it might be intractable by any technique at our disposal. Here again, I slightly alter some notations and correct a confusion.
\medskip
{\it Finally I state an old problem of mine which is probably very difficult and which seems to be unattackable by the methods of probabilistic number theory: denote by $P^+(n)$ the greatest prime factor of $n$. Is it true that the density of integers $n$ satisfying $P^+(n+1)>P^+(n)$ is $\dm$? Is it true that the density of integers for which
$$P^+(n+1)>P^+(n)n^\alpha\eqdef{PP}$$
exists for every $\alpha$? Pomerance and I proved \citer{EP78}  that the upper density of the integers satisfying
$$n^{-\varepsilon}<P^+(n+1)/P^+(n)<n^{\varepsilon}$$
tends to 0 with $\varepsilon$.}
\medskip
Let $E:=\{n\geqslant 1:P^+(n)>P^+(n+1)\}$. The conjecture that $E$
has asymptotic density $\dm$ stems for the general hypothesis that $n$ and $n+1$ should be multiplicatively independent. It lies in the same class of problems than the famous $abc$-conjecture. \par  A general theorem of Hildebrand \citer{Hi85}
implies that
$E$ has positive lower asymptotic density, but I did not check the
numerical value that can be derived from this result.
In \citer{EP78} it is shown that if $N$ is large, then for at least
$0.0099N$ values of $n\le N$ we have $P^+(n)>P^+(n+1)$, and for at least
$0.0099N$ values of $n\le N$ we have $P^+(n)<P^+(n+1)$.  It follows from theorem 1.2 of \citer{BPT05}
 that each inequality occurs on a set of integers $n$ of
lower asymptotic density
$$
\log\left({1\over1-c}\right)-2\int_0^c\log \Big({1-v\over 
1-v-2c}\Big){\d v\over
1-v}
$$
provided $0<c<1/5$.  The maximum of this expression is
greater than $0.05544$, which improves the result from \citer{EP78}.
\par 
In \citer{EP78} it is shown that
$P^+(n)<P^+(n+1)<P^+(n+2)$ holds infinitely often, and it is conjectured
that too $P^+(n)>P^+(n+1)>P^+(n+2)$ holds infinitely often.  This conjecture
was proved by Balog \citer{B01}.
\medskip
Among several, two further very interesting problems are described in Erd\H os' seminal article. I chose not to discuss them in detail since they lie somewhat aside of the main stream of the paper, concentrated on the distribution of divisors and typical multiplicative structure of integers.\par 
 Thus, I only mention (too) briefly the questions of the number $\Phi(x)$ of distinct values of Euler's totient $\varphi(n)$ in $[1,x]$ and that of an infinite sequence $\{p_j\}_{j=1}^{\infty}$ of primes such that $p_{j+1}\md1{p_j}$ $(j\geqslant 1)$.\par 
On the first problem,  a crucial and impressive progress was made by Ford \citer{Fo98}. Improving on results by Pillai \citer{Pi29}, Erd\H os \citer{Er35}, \citer{Er45}, Erd\H os--Hall \citer{EH73}, \citer{EH76}, Pomerance \citer{Po86} and Maier-Pomerance \citer{MP88}, he could show that, for large $x$, we have
$$\Phi(x)\asymp {x\over \log x}\e^{C(\log_3x-\log_4x)^2}(\log_2x)^D(\log_3x)^E,$$
where $C$ and $D$ are positive, explicitly defined constants and $E=D-2C+\dm$.
\par 
On the second problem, Erd\H os asks whether we necessarily
 have $\lim p_j^{1/j}=\infty$ and expresses the belief that $p_j<\exp\{{j(\log j)^{1+o(1)}}\}$ is possible. To my knowledge, both questions are still open. However, Ford, Konyagin and Luca made significant progress in~\citer{FKL10}.\goodbreak 
\bigskip
In conclusion and in the spirit described in the introduction of this article, I hope that this paper will meet two goals. The first is, as for any survey paper, to set records straight, isolate problems and stimulate further research. \par 
Intimately linked to the personality of this so special and so moving (in every sense) man Paul Erd\H os was, the second goal consists in modestly helping to maintain a fair picture of his offering to mathematics. His problems have too often been considered as tricky, disconnected questions. All those who worked with him for some time will agree that, even unformulated, he had in mind the bases of many theories and of even more links between these theories. Now that he can read in the Great Book all answers to his innumerable questions, and indeed select the most elegant ones, no doubt he grins once in a while, realizing how close he has been and pondering how many clues he left for us, even if we still cannot understand them all.
\medskip
\noi{\it Acknowledgements.} The author takes pleasure in expressing here warm thanks to R.~Balasubramanian, N. Bingham, R. de la Bretche, C. Dartyge, I.Z. Ruzsa and T. Stoll for their help during the preparation of this paper.
\goodbreak
\bigskip\medskip
\centerline{\twelvebf References}
\medskip
{\eightpoint\leftskip9mm
\bibtem{AK95} R. Ahlswede \& L.H. Khachatrian, 
Density inequalities for sets of multiples,
{\it J. Number Theory \bf55} \numero2 (1995), 170--180.\par 
\bibtem{B01} A. Balog,  On triplets with descending largest prime
factors, {\it Studia Sci.\ Math.\ Hungar. \bf38} (2001), 45--50.
\bibtem{Be35} F.A. Behrend, Three reviews; of papers by Chowla, Davenport and Erd\H os, {\it Jahrbuch Ÿber die Forschritte der Mathematik \bf 60} (1935), 146--149.\par 
\bibtem{Be48} F.A. Behrend, Generalizations of an inequality of Heilbronn and Rohrbach, {\it Bull.
Amer. Math. Soc.} {\bf 54} (1948), 681--684.\par 
\bibtem{Bes34} A.S. Besicovitch, On the density of certain sequences, {\it Math. Annalen \bf 110} (1934), 336--341.\par 
\bibtem{BPT05} R. de la Bretche, C. Pomerance \& G. Tenenbaum, On products of ratios of consecutive integers, {\it Ramanujan J. \bf9} (2005), 131--138.\par 
\bibtem{BT02} R. de la Bretche \& G. Tenenbaum,  Sur les lois locales de la r\'epartition du
$k$-i\`eme diviseur d'un entier, \PLMS\ (3) \bf84 \rm (2002), 289--323.\par 
\bibtem{BT12} R. de la Bretche \& G. Tenenbaum, Oscillations localisŽes sur les diviseurs, {\it J. London Math. Soc. \rm(2) \bf 85 \rm (2012)} 669Ð693.\par 
\bibtem{BT12a} R. de la Bretche \& G.
Tenenbaum, Moyennes de fonctions arithmŽtiques de formes binaires, {\it Mathematika}, to appear.\par
\bibtem{BT12+} R. de la Bretche \& G.
Tenenbaum, Conjecture de Manin  pour certaines
surfaces de Ch‰telet, preprint.\par 
\bibtem{Ch34} S. Chowla, On abundant numbers, {\it J. Indian Math. Soc.} (2) {\bf 1} (1934), 41--44.\par 
\bibtem{Da33} H. Davenport, †ber numeri abundantes, {\it Sitzungsber. Preu§. Akad. Wiss.}, Phys.-Math. Kl. {\bf 26--29} (1933), 830--837.\par  
\bibtem{DE37} H. Davenport \& P. Erd\H os, On sequences of positive integers, {\it Acta Arith. \bf 2} (1937), 147--151.\par 
\bibtem{DE51} H. Davenport \& P. Erd\H os, On sequences of positive integers, {\it J. Indian Math. Soc. \bf 15} (1951), 19--24. \par 
\bibtem{DKT02} J.-M. De Koninck \& G. Tenenbaum, Sur la loi de r\'epartition du $k$-i\`eme facteur premier d'un entier,
{\it Math. Proc. Camb. Phil. Soc. \bf133} (2002), 191--204. \par                           
\bibtem{Er34} P. Erd\H os, On the density of the abundant numbers, {\it J. London Math. Soc. \bf9}, \numero4 (1934), 278--282.\par 
\bibtem{Er35} P. Erd\H os, On the normal number of prime factors of $p-1$ and some related problems concerning EulerÕs
$\varphi$-function, {\it Quart. J. Math. (Oxford) \bf 6} (1935), 205--213.\par 
\bibtem{Er36} P. Erd\H os, A generalization of a theorem of Besicovitch, \JLMS\ \bf
11 \rm(1936), 92--98.\par 
\bibtem{Er45}  P. Erd\H os, Some remarks on EulerÕs $\varphi$-function and some related problems, {\it Bull. Amer. Math. Soc. \bf51} (1945), 540Ð544.\par 
\bibtem{Er48} P. Erd\H os, On the density of some sequences of integers, {\it Bull. Amer. Math Soc. \bf54} (1948), 685--692.\par 
\bibtem{Er64} P. Erd\H os, On some applications of probability to analysis and number theory, {\it J. London Math. Soc. \bf39} (1964), 692--696.\par 
\bibtem{Er69} P. Erd\H os, On the distribution of prime divisors, {\it Aequationes Math. \bf 2},
(1969), 177--183. \par 
\bibtem{Er74} P. Erd\H os, Problem 218 and solution, {\it Canad. Math. Bull. \bf17} (1974), 621--622.\par 
\bibtem{Er79} P. Erd\H os, Some unconventional problems in number theory, {\it AstŽrique \bf 61} (1979), 73--82, Soc. math. France.\par 
\bibtem{EH73} P. Erd\H os \& R.R. Hall, On the values of EulerÕs $\varphi$-function, {\it Acta Arith. \bf22} (1973), 201--206.\par 
\bibtem{EH76} P. Erd\H os and R.R. Hall, Distinct values of EulerÕs $\varphi$-function, {\it Mathematika \bf23} (1976), 1--3.\par 
\bibtem{EH79} P. Erd\H os \& R.R. Hall, 
The propinquity of divisors, {\it Bull. Lodon Math. Soc. \bf 11} (1979), 304--307.\par 
\bibtem{EK40} P. Erd\H os \& M. Kac,
 The Gaussian law of errors in the theory of additive number theoretic functions,
{\it Amer. J. Math. \bf 62} (1940), 738--742.\par 
\bibtem{EHT94} P. Erd\H os,
R.R. Hall \& G. Tenenbaum, On the densities of sets of multiples, {\it J. reine angew. Math. \bf 454} (1994), 119--141.\par
\bibtem{EN75} P. Erd\H os \& J.-L. Nicolas, RŽpartition des nombres superabondants, {\it Bull. Soc. Math. France \bf103} (1975), 65--90.\par 
\bibtem{EN76} P. Erd\H os \& J.-L. Nicolas, MŽthodes probabilistes et combinatoires en thŽorie des nombres, {\it Bull. sc. math. } (2) {\bf100} (1976), 301--320.\par 
\bibtem{EP78} P. Erd\H os \& C. Pomerance, 
On the largest prime factors of $n$ and $n+1$,
{\it Aequationes Math. \bf17} (1978), no. 2-3, 311--321. \par 
\bibtem{ET81} P. Erd\H os \& G. Tenenbaum, Sur la structure de la suite des diviseurs d'un entier, {\it Ann.
Inst. Fourier (Grenoble) \bf 31}, 1 (1981), 17--37.\par
\bibtem{ET83} P. Erd\H os \& G. Tenenbaum, Sur les diviseurs cons\'ecutifs d'un entier, {\it Bull. Soc. Math. de
France \bf 111} (1983), 125--145. \par 
\bibtem{ET89a} P. Erd\H os \& G. Tenenbaum, Sur les densit\'es de certaines suites d'entiers, {\it
Proc. London Math. Soc.,} {\rm (3)} {\bf 59} (1989), 417--438.
\bibtem{ET89} P. Erd\H os \& G. Tenenbaum, Sur les fonctions arithm\'etiques
li\'ees aux diviseurs cons\'ecutifs, {\it J. Number Theory, \bf 31} (1989), 285--311.\par 
\bibtem{Fo98} K. Ford, 
The distribution of totients, {\it
Ramanujan J. \bf2} (1998), \numeros 1-2, 67--151. \par 
\bibtem{Fo08} K. Ford, 
The distribution of integers with a divisor in a given interval. {\it Ann. of Math. (2) \bf168}, \numero 2 (2008), 367Ð433. \par 
\bibtem{FKL10} K. Ford, S.V. Konyagin \& F. Luca, Prime chains and Pratt trees, {\it Geom. Funct. Anal. \bf20} \numero 5 (2010), 1231--1258.\par 
\bibtem{HR74} H. Halberstam \& H.-E. Richert,
 {\it Sieve Methods}, Academic Press, London, New York, San Francisco, 1974.\par 
\bibtem{Ha90} R.R. Hall, Sets of multiples and Behrend sequences, {\it A tribute to Paul Erd\H os}, 249--258, Cambridge Univ. Press, Cambridge, 1990.\par 
\bibtem{Ha96} R.R. Hall, {\it Sets of multiples},
Cambridge Tracts in Mathematics, 118, Cambridge University Press, Cambridge, 1996.\par 
\bibtem{Ha12} R.R. Hall, Private communication, November 11, 2012.\par 
\bibtem{HT82} R.R. Hall \& G. Tenenbaum, On the average and normal orders of Hooley's $\Delta$-function,
{\it J. London Math. Soc} (2) {\bf 25} (1982), 392--406.\par 
\bibtem{HT86} R.R. Hall \& G. Tenenbaum, Les ensembles de multiples et la densit\'e divisorielle, {\it J.
Number Theory \bf 22} (1986), 308--333.\par 
\bibtem{HT88} R.R. Hall \& G. Tenenbaum, {\it Divisors}, Cambridge tracts in
mathematics 90, Cambridge University Press (1988, paperback ed. 2008).\par 
\bibtem{HT92} R.R. Hall \& G. Tenenbaum, On Behrend sequences, 
{\it Math. Proc. Camb. Phil. Soc. \bf 112} (1992), 467--482.\par 
\bibtem{He12} K. Henriot, Nair--Tenenbaum bounds
uniform with respect to the discriminant, {\it Math. Proc. Camb. Phil. Soc. \bf152} (2012),\numero3, 405--424.\par 
\bibtem{Hi85} A. Hildebrand, On a conjecture of Balog,
{\it Proc. Amer. Math. Soc. \bf95} , n$^\circ\thinspace$4 (1985), 517--523.\par 
\bibtem{Ho79} C. Hooley, A new technique and its applications to the theory of numbers, {\it Proc. London Math. Soc.} (3) {\bf 38} (1979), 115--151.\par
\bibtem{KT04} S. Kerner \& G. Tenenbaum, Sur la r\'epartition divisorielle normale de $\theta d\,(\mod1)$, {\it Math. Proc. Camb. Phil. Soc. \bf137} (2004), 255--272.\par 
\bibtem{Le37} P. LŽvy, {\it Th\'eorie de l'addition des variables al\'eatoires}, 
Gauthier-Villars, Paris (1937, 2nd ed. 1954). \par 
\bibtem{Ma87} H. Maier, On the Mšbius function,
{\it Trans. Amer. Math. Soc. \bf 301},
\numero 2 (1987), 649--664.\par
\bibtem{MP88} H. Maier \& C. Pomerance, On the number of distinct values of EulerÕs $\varphi$-function, {\it Acta Arith. \bf49} (1988), 263--275.\par 
\bibtem{MT84} H. Maier \& G. Tenenbaum, On the set of divisors of an integer, {\it Invent. Math. \bf 76} (1984), 121--128.\par 
\bibtem{MT11} H. Maier \& G. Tenenbaum, On the normal concentration of divisors, 2, {\it Math. Proc. Camb. Phil. Soc. \bf147} \numero3 (2009), 593--614.\par
\bibtem{NT98} M. Nair \& G. Tenenbaum, Short sums
of certain arithmetic functions, {\it Acta Math.}
{\bf 180}, (1998), 119--144.\par 
\bibtem{Pi29} S. Pillai, On some functions connected with $\varphi(n)$, {\it Bull. Amer. Math. Soc. \bf35} (1929), 832--836.\par 
\bibtem{PS98} G. P—lya \& G. Szeg\H o,  {\it Problems and theorems in analysis, II}, Classics in Mathematics, Springer-Verlag, Berlin, 1998.\par 
\bibtem{Po86} C. Pomerance, On the distribution of the values of EulerÕs function, {\it Acta Arith. \bf47} (1986), 63--70.\par 
\bibtem{Ra95} A. Raouj, 
Sur la densitŽ de certains ensembles de multiples, I, II, 
{\it Acta Arith. \bf69}, \numero2 (1995), 121Ð152, 171Ð188.\par  
\bibtem{RST11} A. Raouj, A. Stef \& G. Tenenbaum, Mesures quadratiques de la proximit\'e des diviseurs, {\it Math. Proc. Camb. Phil. Soc. \bf150} (2011), 73--96.\par 
\bibtem{Ro11} O. Robert, 
Sur le nombre des entiers reprŽsentables comme somme de trois puissances, {\it
Acta Arith. \bf149} \numero1 (2011), 1Ð21. \par 
\bibtem{RT96} I.Z. Ruzsa \& G. Tenenbaum, A note on Behrend sequences, {\it Acta Math. Hung.} 
\bf 72\rm, \numero 4 (1996), 327--337.\par 
\bibtem{St92} A. Stef, {\it L'ensemble exceptionnel dans la conjecture d'Erd\H os concernant la proximitŽ des diviseurs}, Thse de doctorat de l'UniversitŽ Nancy 1, UFR STMIA, juin 1992.\par 
\bibtem{Te82} G. Tenenbaum, Sur la densit\'e divisorielle d'une suite d'entiers, {\it J. Number Theory \bf
15}, n$^{\circ}$ 3 (1982), 331-346.\par 
\bibtem{Te84} G. Tenenbaum, Sur la probabilit\'e qu'un entier poss\`ede un diviseur dans un intervalle donn\'e,      
{\it Compositio Math. \bf 51} (1984), 243--263.\par 
\bibtem{Te85} G. Tenenbaum, Sur la concentration moyenne des diviseurs,  {\it Comment. Math. Helvetici \bf 60}
(1985), 411--428.
\par
\bibtem{Te86} G. Tenenbaum, Fonctions $\Delta$ de Hooley et applications, {\it S\'eminaire de Th\'eorie des
nombres, Paris 1984-85, Prog. Math. \bf 63} (1986), 225--239.\par 
\bibtem{Te87} G. Tenenbaum, Un probl\`eme de probabilit\'e conditionnelle en Arithm\'etique, {\it Acta
Arith. \bf 49} (1987), 165--187.\par 
\bibtem{Te90} G. Tenenbaum, Sur une question d'Erd\H os et Schinzel, II, {\it Inventiones Math.} {\bf 99}
(1990), 215--224.\par 
\bibtem{Te96} G. Tenenbaum, On block Behrend sequences, {\it Math. Proc. Camb. Phil. Soc. \bf 120} (1996),
355--367.\par 
\bibtem{Te96a} G. Tenenbaum, Uniform distribution on divisors and Behrend sequences, {\it L'Enseignement
Math\'e\-matique \bf 42} (1996), 153--197.\par 
\bibtem{Te08} G. Tenenbaum, {\it Introduction \`a la th\'eorie analytique et probabiliste des
nombres}, third ed., coll. ƒchelles, Belin, 2008, 592 \pp.\par 
\bibtem{RCV85} R.C. Vaughan, 
Sur le problme de Waring pour les cubes, {\it
C. R. Acad. Sci. Paris SŽr. I Math. \bf 301}, \numero 6 (1985), 253Ð255.\par }
\vskip 5mm
{\sevenrm\baselineskip9pt
G\'erald Tenenbaum\par
Institut \'Elie Cartan\par 
Universit\'e Henri Poincar\'e--Nancy 1\par
 BP 239\par
54506 Vand\oe uvre Cedex\par
 France
\smallskip
internet: \seventt gerald.tenenbaum@iecn.u-nancy.fr\par}

\bye